\documentclass[a4paper,11pt]{article}%
\usepackage{amsmath}
\usepackage{amsfonts}
\usepackage{amssymb}
\usepackage{graphicx}%
\newtheorem{theorem}{Theorem}
\newtheorem{proposition}{Proposition}
\newtheorem{remark}{Remark}
\newtheorem{corollary}{Corollary}
\newtheorem{definition}{Definition}

\newtheorem{lemma}{Lemma}
\newenvironment{proof}[1][Proof]{\noindent\textbf{#1.} }{\ \rule{0.5em}{0.5em}}
\begin{document}

\title{Representation of small ball probabilities in Hilbert space and lower bound in
regression for functional data.}
\author{Andr\'{e} Mas}
\maketitle

\begin{abstract}
Let $S=\sum_{i=1}^{+\infty}\lambda_{i}Z_{i}$ where the $Z_{i}$'s are i.d.d.
positive with $\mathbb{E}\left\vert Z\right\vert ^{3}<+\infty$ and $\left(
\lambda_{i}\right)  _{i\in\mathbb{N}}$ a positive nonincreasing sequence such
that $\sum\lambda_{i}<+\infty$. We study the small ball probability
$\mathbb{P}\left(  S<\varepsilon\right)  $ when $\varepsilon\downarrow0$. We
start from a result by Lifshits (1997) who computed this probability by means
of the Laplace transform of $S$. We prove that $\mathbb{P}\left(
S<\cdot\right)  $ belongs to a class of functions introduced by de Haan,
well-known in extreme value theory, the class of Gamma-varying functions, for
which an exponential-integral representation is available. This approach
allows to derive bounds for the rate in nonparametric regression for
functional data at a fixed point $x_{0}$ : $\mathbb{E}\left(  y|X=x_{0}%
\right)  $ where $\left(  y_{i},X_{i}\right)  _{1\leq i\leq n}$ is a sample in
$\left(  \mathbb{R},\mathcal{F}\right)  $ and $\mathcal{F}$ is some space of
functions. It turns out that, in a general framework, the minimax lower bound
for the risk is of order $\left(  \log n\right)  ^{-\tau}$ for some $\tau>0$
depending on the regularity of the data and polynomial rates cannot be achieved.

\end{abstract}


\section{Preliminaries}

The three following susbections are independent. The first gives some basic
material about small ball probability. The second collects classical results
from extreme value theory as well as the definition of the class $\Gamma_{0}$
which is then briefly described. The third introduces the nonparametric
regression model for functional data and simply raises the problems attached
to obtaining sharp bounds for the quadratic risk at a fixed point. The notions
encountered in this long introduction though intially distinct from each other
merge in the sequel of this work and give birth to the main results. Some
proofs are given in the last section.

\subsection{About non-shifted and shifted small ball problems}

Small ball problems could generally be stated the following way : consider a
random variable $X$ with values in a general normed space $\left(
E,\left\Vert \cdot\right\Vert \right)  $ (which may not be finite-dimensional)
and estimate $\mathbb{P}\left(  \left\Vert X\right\Vert <\varepsilon\right)  $
for small values of $\varepsilon$. This issue may be viewed as a counterpart
of the large deviations or concentration problems (where $\mathbb{P}\left(
\left\Vert X\right\Vert >M\right)  $ is studied for large $M$) and the terms
"small deviations" or "lower tail behaviour" are sometimes encountered to name
small ball problems. The core of the literature on small ball problems focuses
on Gaussian random variables. The survey by Li and Shao (2001) is a complete
state of the art, introducing the main concepts and providing numerous
references. Another reference is Chapter 18 of Lifshits (1995) entirely
devoted to Gaussian random functions. Much attention has been given to
Brownian motion (when $\left(  E,\left\Vert \cdot\right\Vert \right)  =\left(
C\left(  0.1\right)  ,\left\vert \cdot\right\vert _{\infty}\right)  $) or its
relatives (fractional Bronwian motion, Bronwian sheet, etc). The case of
stable random elements was also investigated (see for instance Li, Linde
(2004), Aurzada, Lifshits, Linde (2009)). Another issue is related to the
norm. Indeed in infinite dimensional spaces, norms or metrics are not
equivalent and this may influence the local behaviour of $\mathbb{P}\left(
\left\Vert X\right\Vert <\varepsilon\right)  $.

A more general question could be the shifted small ball probability
$\mathbb{P}\left(  \left\Vert X-x_{0}\right\Vert <\varepsilon\right)  $ for a
fixed $x_{0}$. A stumbling stone arises from the shift $x_{0}$. It turns out
that, in general, computations cannot be carried out for any $x_{0}$. Several
works focus on expliciting the set of those $x_{0}$ for which the shifted
small ball probability may be computed from the non-shifted one (when
$x_{0}=0$). We refer to Borell (1976) or Li and Linde (1993) for instance. A
classical example stems from the situation where $\mathbb{P}_{X-x_{0}}%
\ll\mathbb{P}_{X}$ where $\mathbb{P}_{X}$ denotes the probability distribution
induced by the random element $X$. The classical Cameron-Martin's theorem for
Brownian motion illustrates this case for instance. Absolute regularity yields
:%
\[
\mathbb{P}\left(  \left\Vert X-x_{0}\right\Vert <\varepsilon\right)
=\int_{B\left(  0,\varepsilon\right)  }\mathbb{P}_{X-x_{0}}\left(  dx\right)
=\int_{B\left(  0,\varepsilon\right)  }f_{x_{0}}\left(  x\right)
\mathbb{P}_{X}\left(  dx\right)
\]
where $f_{x_{0}}=d\mathbb{P}_{X-x_{0}}/d\mathbb{P}_{X}$ and $B\left(
0,\varepsilon\right)  $ stands for the ball centered at $0$ with radius
$\varepsilon$. When $f_{x_{0}}$ is regular enough in a neighborhood of zero :%
\begin{equation}
F_{x_{0}}\left(  h\right)  =\int_{B\left(  0,\varepsilon\right)  }%
\mathbb{P}_{X-x_{0}}\left(  dx\right)  =\int_{B\left(  0,\varepsilon\right)
}f_{x_{0}}\left(  x\right)  \mathbb{P}_{X}\left(  dx\right)  \sim f_{x_{0}%
}\left(  0\right)  F\left(  h\right)  \label{eqshsbp}%
\end{equation}

About this fact see Proposition 2.1 in de Acosta (1983). In general the
sharpness of existing results may vary, depending on the triplet $\left(
\left(  E,\left\Vert \cdot\right\Vert \right)  ,\mathbb{P}_{X},x_{0}\right)  $
under consideration. In fact there are only few spaces for which the local
behaviour of $\mathbb{P}\left(  \left\Vert X-x_{0}\right\Vert <\varepsilon
\right)  $ is explicitely described. Quite often lower and upper bounds are
computed so that :%
\[
\mathbb{P}\left(  \left\Vert X-x_{0}\right\Vert <\varepsilon\right)
\asymp\varphi_{x_{0}}\left(  \varepsilon\right)
\]
where $\varphi_{x_{0}}$ is known and $f\asymp g$ means here that for some
positive constants $c^{-}$ and $c^{+}$ the positive functions $f$ and $g$
satisfy :%
\[
0<c^{-}\leq\liminf_{0}\frac{f}{g}\leq\limsup_{0}\frac{f}{g}\leq c^{+}.
\]
Sometimes only one of these bounds is accessible or needed.

It is worth noting or recalling a few crucial features of small deviations
techniques. The Laplace transform, as well as in large deviations problems, is
a major tool when coupled with the saddlepoint method. Small deviations are
intimately connected with the entropy of the unit ball of the reproducing
kernel Hilbert space associated with $X$, with the $l$-approximation numbers
of $X$ (i.e. the rate of approximation of $X$ by a finite dimensional random
variable, see Li, Linde (1999)) or to the degree of compactness of linear
operators generating $X$ (see Li, Linde (2004)). All these notions are clearly
connected to the regularity of the process $X,$ when $X$ is a process.

Applications of small ball probabilities are numerous : they appear when
studying rates of convergence in the Law of the Iterated Logarithm (see
Talagrand (1992), Kuelbs, Li, Linde (1994)) or the rate of escape of the
Brownian motion (see Erickson (1980)). They even surprisingly provide a
sufficient condition for the CLT (see Ledoux, Talagrand (1991), Theorem 10.13
p.289). However small ball problems remained until nowadays a matter
essentially reserved to probability theory. However Van der Vaart and van
Zanten (2007 and 2008) found applications of small ball techniques to Bayesian
statistics It turns out that this topic may be also of interest in another
area of statistics : functional data analysis. FDA for short extends classical
statistical models designed for vectors to the situation when the data are
functions or curves. One of the concern may be summarized this way : since
Lebesgue's density of an infinite-dimensional random $X$ does not exist, all
the inference techniques based on the density cannot hold anymore. In this
framework, the small ball probabilities appear as a natural counterpart and
should be investigated with much care. We illustrate this fact by pointing out
an elementary example -kernel methods- in the next subsection below.

First let us precise the $l_{2}$ framework. Consider $X$ a random variable
defined the following way :%
\begin{equation}
X=\left(  \sqrt{\lambda_{1}}x_{1},\sqrt{\lambda_{2}}x_{2},...\right)
\label{lou}%
\end{equation}
where $\left(  \lambda_{i}\right)  _{1\leq i\leq n}$ is a real positive
sequence arranged in a non-decreasing order such that $\sum_{i=1}^{+\infty
}\lambda_{i}<+\infty$ and $\left(  x_{i}\right)  _{1\leq i\leq n}$ is sequence
of real independent and identically distributed random variables with null
expectation. From Kolmogorov's $0-1$ law it is straightforward to see that $X$
exists as a $l_{2}$-valued random element. The square norm of $X$ is
$S=\sum_{i=1}^{+\infty}\lambda_{i}x_{i}^{2}.$

The small ball problem consists here in estimating for different choices of
the sequence $\left(  \lambda_{i}\right)  _{i\in\mathbb{N}}$ and $\left(
x_{i}\right)  _{i\in\mathbb{N}}$ the probability $\mathbb{P}\left(
S<r\right)  $ when $r$ tends to zero. The latter probability is expected to
depend on the $\lambda_{i}$'s. About this fact we refer to Dunker, Lifshits,
Linde (1998).

The inspection of the case $E=l_{2}$ is motivated by the application to
functional statistics mentioned in the paragraph above. Indeed random
functions are often reconstructed by interpolation techniques, like splines or
wavelets, in Hilbert spaces such as $L^{2}\left(  \left[  0,T\right]  \right)
$ or the Sobolev space $W^{m,2}\left(  \left[  0,T\right]  \right)  ,$
$m\in\mathbb{N}$. Then the random element $X$ is valued in a separable Hilbert
space $\mathcal{H}$ and all these Hilbert spaces of functions are
isometrically isomorphic to $l_{2}.$ In this framework a useful tool is the
so-called Karhunen-Lo\`{e}ve decomposition (sometimes refered to as Principal
Orthogonal Decomposition in other area of mathematics such as PDEs). Any
centered random function $X$ will be represented by its coordinates in a basis
of eigenvectors of the covariance operator $\mathcal{C}_{X}=\mathbb{E}\left[
X\otimes X\right]  $. When $e_{i}$'s are the eigenvectors of $\mathcal{C}_{X}$
and $\lambda_{i}$ the associated eigenvalues
\begin{equation}
X=\sum_{i=1}^{+\infty}\sqrt{\lambda_{i}}x_{i}e_{i} \label{KL}%
\end{equation}
where the $x_{i}$'s are uncorrelated real random variables. The $x_{i}$'s are
actually always independent when $X$ is Gaussian and are assumed to be in most
settings. The $l_{2}$ random element defined in (\ref{lou}) is formally
identifiable with this Karhunen-Lo\`{e}ve decomposition familiar in Functional
Data Analysis.

Historically the description of the exact behaviour of Gaussian small ball
probability in Hilbert space is due to Sytaya (1974). However we borrow the
notations from Lifshits (1997) who extendend Sytaya's results in several
directions amongts which the non-Gaussian framework. First in order to
alleviate notations set once and for all :%
\begin{equation}
S=\sum_{i=1}^{+\infty}\lambda_{i}Z_{i} \label{sbpl2}%
\end{equation}
where $\lambda_{i}>0$ are arranged in decreasing order with $\sum
_{i=1}^{+\infty}\lambda_{i}<+\infty$ and $Z_{i}$ are positive random variables
(they stand for the $x_{i}^{2}$'s above). For the sake of completeness and
since the main theorems of this work heavily rely on his results we recall
them. In the previously mentioned article Lifshits proved that :%
\begin{equation}
\mathbb{P}\left(  S<r\right)  \underset{r\rightarrow0}{\sim}\frac{1}%
{\sqrt{2\pi}}\frac{1}{\gamma\sigma}\exp\left(  \gamma r\right)  \Lambda\left(
\gamma\right)  \label{Dem}%
\end{equation}
where $\gamma$ and $\sigma$ are functions of $r$ defined below and
$\Lambda\left(  \gamma\right)  =\mathbb{E}\exp\left(  -\gamma S\right)  $ is
the Laplace transform of $S$ evaluated at $\gamma\left(  r\right)  $. The
definitions of $\gamma$ and $\sigma$ are implicit. Let $S_{\gamma}$ be the
Esscher transform of $S$ that is the random variable with distribution
$\exp\left(  -\gamma x\right)  \mathbb{P}_{S}\left(  dx\right)  /\Lambda
\left(  \gamma\right)  $. Then set :%
\begin{align}
r  &  =\mathbb{E}\left[  S_{\gamma}\right]  =-\frac{\partial\log\Lambda\left(
\gamma\right)  }{\partial\gamma},\label{mu}\\
\sigma^{2}  &  =\mathbb{V}\left[  S_{\gamma}\right]  =\frac{\partial^{2}%
\log\Lambda\left(  \gamma\right)  }{\partial\gamma^{2}}. \label{D4}%
\end{align}
where $\mathbb{V}$ denotes variance. Without further assumption on the
$\lambda_{i}$'s $\mathbb{P}\left(  S<r\right)  $ cannot be made more explicit.
This is done for instance in Dunker, Lifshits, Linde (1998) where these author
considered the case of $\lambda_{i}$ with polynomial and exponential decay.
Due to the remark below (\ref{KL}) we will sometimes refer to $\mathbb{P}%
\left(  S<r\right)  $ as a small ball probability for an $l_{2}$-valued random element.

The article is organized as follows. The next subsection develops some aspects
of mathematical statistics which motivate this approach on small ball
problems. Then a class of functions which appears in extremes value theory
-the class $\Gamma_{0}$- is introduced in the next section. Our main theorem
shows that small ball probabilities of $l_{2}$ random elements (hence of
random functions belonging to a Hilbert space) belong to the class $\Gamma
_{0}$. We then show how this result may be used to solving the statistical
issues mentioned earlier. In particular we prove the the optimal rate of
convergence in nonparametic regression for functional variables is always
slower than any power of $n$. The derivations of the main results are
collected in the last part of the article.

\subsection{The class $\Gamma_{0}$}

The theory of extremes is another well-known topic connecting probability
theory, mathematical statistics and real analysis through regular variation
and Karamata's theory. The foundations of extreme value theory may be
illustrated by the famous Fisher-Tippett theorem (see Fisher, Tippett (1928)
and Gnedenko (1943)). This classical result assesses that whenever
$U_{1},...,U_{n}$ is an i.id. sample of real random variables, $M_{n}%
=\max\left\{  U_{1},...,U_{n}\right\}  $ belongs to the domain of attraction
of $G$, where $G$ has same type as one of the three distributions Gumbel,
Frechet and Weibull. The Gumbel law, also named double exponential
distribution, with cumulative distribution function $\Lambda\left(  x\right)
=\exp\left(  -\exp\left(  -x\right)  \right)  $ defines the so-called "domain
of attraction of the third type". Laurens de Haan (1971) characterized the
(cumulative) distribution functions of $U$ such that $M_{n}$ belongs to the
domain of attraction of $\Lambda$. We give this result below.

\textbf{Theorem (de Haan, 1971) :} \textit{If }$F$\textit{ is the cumulative
distribution function of a real random variable }$X$\textit{ which belongs to
the domain of attraction of the third type (Gumbel) there exists a measurable
function }$\rho:\mathbb{R}\rightarrow\mathbb{R}^{+}$\textit{, called the
auxiliary function of }$F$, \textit{such that :}%
\[
\lim_{s\uparrow x_{+}}\frac{\overline{F}\left(  s+x\rho\left(  s\right)
\right)  }{\overline{F}\left(  s\right)  }=\exp\left(  -x\right)
\]
\textit{where }$\overline{F}\left(  s\right)  =1-F\left(  s\right)  ,$\textit{
}$x_{+}=\sup\left\{  x:F\left(  x\right)  <1\right\}  $.

This property was intially introduced by de Haan as a "Form of Regular
Variation" (see the title of his article). This class of distribution function
is referred to as de Haan's Gamma class in the book by Bingham, Goldie and
Teugels (1987) and within this article. In the latter book the definition is
slightly different from the one given above. Gamma-variation is defined at
infinity and for non-decreasing functions which comes down to taking
$x_{+}=+\infty$ and taking $\exp\left(  x\right)  $ instead of $\exp\left(
-x\right)  $ in the display above. Surprisingly, in their book as well as in
de Haan's article no examples of functions belonging to $\Gamma$ is given. The
cumulative distribution function function of the Gaussian distribution belongs
to this class with $x_{+}=+\infty$ and $\rho\left(  s\right)  =1/s$.

Since we focus on the local behaviour at zero of the cumulative distribution
function function of a real valued random variable we have to modifiy again
slightly the definitions above. We introduce the class $\Gamma_{0}$ and
feature some of its properties below. We borrow most of our notations from
Bingham, Goldie and Teugels (1987) which differ from those of de Haan.

\begin{definition}
\label{def}The class $\Gamma_{0}$ consists of those functions $F:\mathbb{R}%
\rightarrow\mathbb{R}^{+}$ null over $\left(  -\infty,0\right]  $, non
decreasing with $F\left(  0\right)  =0$ and right-continuous for which there
exists a continuous non decreasing function $\rho:\mathcal{V}^{+}%
\rightarrow\mathbb{R}^{+}$, defined on some a right-neighborhood of zero
$\mathcal{V}^{+}$ such that $\rho\left(  0\right)  =0$ and for all
$x\in\mathbb{R},$%
\begin{equation}
\lim_{s\downarrow0^{+}}\frac{F\left(  s+x\rho\left(  s\right)  \right)
}{F\left(  s\right)  }=\exp\left(  x\right)  \label{deflim}%
\end{equation}
The function $\rho$ is called the auxiliary function of $F$.
\end{definition}

The properties of the auxiliary function are crucial.

\begin{proposition}
\label{cabrel}From Definition \ref{def} above we deduce that : $\rho\left(
s\right)  /s\rightarrow0$ as $s\rightarrow0$ and $\rho$ is self-neglecting
which means that :%
\[
\frac{\rho\left(  s+x\rho\left(  s\right)  \right)  }{\rho\left(  s\right)
}\overset{s\rightarrow0}{\rightarrow}1
\]
locally uniformly in $x\in\mathbb{R}$.
\end{proposition}

\begin{remark}
When the property in the proposition above does not hold locally uniformly but
only pointwise the function is called Beurling slowly varying. Assuming that
$\rho$ is continuous in Definition \ref{def} yields local uniformity and
enables to consider a self-neglecting $\rho.$
\end{remark}

The class $\Gamma_{0}$ is subject to an exponential-integral representation.
In fact the following Theorem asserts that the local behaviour at $0$ of any
$F$ in $\Gamma_{0}$ depends only on the auxiliary mapping $\rho$.

\begin{theorem}
\label{repres}Let $F$ belong to $\Gamma_{0}$ with self-neglecting auxiliary
function $\rho$ then when $s\rightarrow0$ :%
\begin{equation}
F\left(  s\right)  =\exp\left\{  \eta\left(  s\right)  -\int_{s}^{1}\frac
{1}{\rho\left(  t\right)  }dt\right\}  \label{zazie}%
\end{equation}
with $\eta\left(  s\right)  \rightarrow c\in\mathbb{R}$ and the auxiliary
function $\rho$ is unique up to asymptotic equivalence and may be taken as
$\int_{0}^{s}F\left(  t\right)  dt/F\left(  s\right)  .$ Besides%
\begin{equation}
F\left(  \lambda s\right)  /F\left(  s\right)  \rightarrow\left\{
\begin{tabular}
[c]{ll}%
$\infty$ & $\left(  \lambda>1\right)  $\\
$1$ & $\left(  \lambda=1\right)  $\\
$0$ & $\left(  \lambda<1\right)  $%
\end{tabular}
\ \ \ \ \right.  \quad\text{as }s\rightarrow0. \label{zaza}%
\end{equation}

\end{theorem}

\begin{remark}
The upper bound $1$ in the integral in display (\ref{zazie}) is unimportant
and may be replaced by any positive number. Then the function $\eta$ will
change as well.
\end{remark}

The proof of Proposition \ref{cabrel} as well as Theorem \ref{repres} are
inspired from the proofs of Lemma 3.10.1, Proposition 3.10.3 and Theorem
3.10.8 in Bingham et al (1987) and will be omitted.\newline Let us also
mention that Ga\"{\i}ffas (2005) proposed to model locally the density of
sparse data by gamma-varying functions. This is another statistical
application for $\Gamma_{0}$. It is simple to construct explicit examples of
functions in $\Gamma_{0}$ by tuning the auxiliary function $\rho$ and taking
$\eta\left(  \cdot\right)  =0$ in (\ref{zazie}). For instance taking $\rho
_{1}\left(  t\right)  =t^{m}$ (with $m>1$) gives $\mathcal{F}_{1}\left(
s\right)  =\exp\left(  -1/s^{m-1}\right)  $. Now taking $\rho_{2}\left(
t\right)  =-t/\log\left(  t\right)  $ yields $\mathcal{F}_{2}\left(  s\right)
=\exp\left(  -\left[  \log\left(  s\right)  \right]  ^{2}\right)  .$ Obviously
constants may be added in front of or within the exponential. The next
Proposition seems to show a specific feature of the class $\Gamma_{0}$.

\begin{proposition}
\label{skye}Let $F$ belong to $\Gamma_{0}$. Then for all integer $p$
$F^{\left(  p\right)  }\left(  0\right)  =0$ where $F^{\left(  p\right)  }$
denotes the derivative of order $p$ of $F$.
\end{proposition}

\subsection{The nonparametric regression model for functional data}

As a last part of this introduction we shift from small ball problems and
extreme theory to statistics for functional data. This recent domain of
statistics has been receiving increasing interest and was boosted by
computational advances. We briefly recall that the main purpose of functional
data analysis (FDA) is to model and study datasets where observations are of
functional nature (usually observed on a grid then smoothed, approximated and
reconstructed by projection on accurate basis ). We refer to the monographs by
Ramsay and Silverman (2005) and Ferraty and Vieu (2006) for an overview of
this topic. Along the past decade some authors turned their attention to the
question of modelizing probability distribution for curve-data with
applications in statistics : Dabo-Niang (2002), Hall and Heckman (2002)
Delaigle and Hall (2010) in a general setting then Dabo-Niang, Ferraty and
Vieu (2004 and 2006), Ferraty, Go\"{\i}a and Vieu (2007) with applications to
classifications through modal curves for instance. Consider the regression
problem with functional data as inputs :%
\begin{equation}
y=r\left(  X\right)  +\varepsilon\label{reg-mod}%
\end{equation}
where $y,\varepsilon$ are real with $\varepsilon$ centered whose variance is
denoted $\sigma_{\varepsilon}^{2}$, $X$ belongs to the Hilbert space
$\mathcal{H}$ and $r$ is a function from $\mathcal{H}$ to $\mathbb{R}$. The
space $\mathcal{H}$ may be chosen to be $L^{2}\left(  T\right)  $ where $T$ is
a compact set in the Euclidean space or som Sobolev space $\mathbf{H}^{2,m}$.
It is endowed with an inner product $\left\langle \cdot,\cdot\right\rangle $
inducing a norm $\left\Vert \cdot\right\Vert .$ Estimating the regression
function at a fixed point $x_{0}$ namely $r\left(  x_{0}\right)
=\mathbb{E}\left(  y|X=x_{0}\right)  $ is possible by a classical
Nadarya-Watson approach (see Tsybakov (2004) for a general presentation in the
finite dimensional setting and Ferraty Vieu (2006) for implementation on
functional data). This model was studied for instance in Ferraty, Vieu (2004)
and asymptotic results were derived in Ferraty, Mas, Vieu (2007) like a first
upper bound for the quadratic risk. It seems that an equivalent of
projection-based estimate in this model has not been introduced yet, certainly
due to a lack of theoretical results on approximation theory for functions
defined on a Hilbert space. The linear regression model $y=\int X\left(
s\right)  \beta s+\varepsilon$ has been extensively investigated in the last
years and several authors proved optimality results like for instance Hall and
Horowitz (2007), Crambes, Kneip, Sarda (2009) or Cardot and Johannes (2010)
(see also references therein these works). It seems that the optimal (in
minimax sense) asymptotic risk has not been obtained yet in the more general
model (\ref{reg-mod}). The behaviour of the small ball probability was a
stumbling stone hard to circumvent.

An adapted Nadaraya-Watson estimate reads :%
\[
\widehat{r}\left(  x_{0}\right)  =\frac{\sum_{i=1}^{n}y_{i}K\left(  \left\Vert
X_{i}-x_{0}\right\Vert /h\right)  }{\sum_{i=1}^{n}K\left(  \left\Vert
X_{i}-x_{0}\right\Vert /h\right)  }%
\]
where $x_{0}$ is a fixed point of the space, $K$ is a kernel, that is a
mesurable, unilateral (defined on $\mathbb{R}^{+}$) positive function with
$\int K=1$ and $h$ is a nonnegative number tending to $0$ (the bandwidth).
Considering the $L^{2}$-risk at a fixed point $x_{0}$ leads to a bias-variance
decomposition :%
\[
\mathcal{R}_{n}\left(  x_{0}\right)  =\mathbb{E}\left(  \widehat{r}\left(
x_{0}\right)  -r\left(  x_{0}\right)  ^{2}\right)  =\mathcal{B}_{n}\left(
x_{0}\right)  +\mathcal{V}_{n}\left(  x_{0}\right)
\]
with%
\begin{align}
\mathcal{B}_{n}\left(  x_{0}\right)   &  =\left\{  \mathbb{E}\frac{\sum
_{i=1}^{n}\left[  r\left(  X_{i}\right)  -r\left(  x_{0}\right)  \right]
K\left(  \left\Vert X_{i}-x_{0}\right\Vert /h\right)  }{\sum_{i=1}^{n}K\left(
\left\Vert X_{i}-x_{0}\right\Vert /h\right)  }\right\}  ^{2}\label{bias}\\
\mathcal{V}_{n}\left(  x_{0}\right)   &  =\mathbb{E}\left[  \frac{\sum
_{i=1}^{n}\varepsilon_{i}K\left(  \left\Vert X_{i}-x_{0}\right\Vert /h\right)
}{\sum_{i=1}^{n}K\left(  \left\Vert X_{i}-x_{0}\right\Vert /h\right)
}\right]  ^{2}\label{var}%
\end{align}

where $\varepsilon_{i}=y_{i}-r\left(  X_{i}\right)  $.

\begin{lemma}
\label{risk}The following holds for the two components of the risk at a fixed
point $x_{0}$ of the kernel estimator $\widehat{r}\left(  x_{0}\right)  $ :%
\begin{align}
\mathcal{V}_{n}\left(  x_{0}\right)   &  \sim\frac{\sigma_{\varepsilon}^{2}%
}{n}\frac{\mathbb{E}K^{2}\left(  \left\Vert X-x_{0}\right\Vert /h\right)
}{\left[  \mathbb{E}K\left(  \left\Vert X-x_{0}\right\Vert /h\right)  \right]
^{2}}\label{v}\\
\mathcal{B}_{n}\left(  x_{0}\right)   &  \sim\frac{1}{\mathbb{E}^{2}K\left(
\left\Vert X-x_{0}\right\Vert /h\right)  }\left(  \sum_{i=1}^{+\infty
}\mathbf{b}_{i}\mathbb{E}\left[  \left\langle X,e_{i}\right\rangle
^{2}K\left(  \left\Vert X\right\Vert /h\right)  \right]  \right)  ^{2}
\label{b}%
\end{align}
where the $\mathbf{b}_{i}$'s are positive and non random constants.
\end{lemma}

The sequence $\left(  \mathbf{b}_{i}\right)  _{i\in\mathbb{N}}$ is not given
here because it depends on several parameters which will be introduced later.
The proof of this lemma will not be explicitely carried out. It will be
encapsulated in the proof of Proposition \ref{risk2} which is more precise
about the bounds (\ref{v}) and (\ref{b}). We keep in mind that the
bias-variance decomposition of the risk is essentially based on the
computation of two sorts of moments : $\mathbb{E}K\left(  \left\Vert
X-x_{0}\right\Vert /h\right)  $ and $\mathbb{E}\left[  \left\langle
X,e_{i}\right\rangle ^{2}K\left(  \left\Vert X\right\Vert /h\right)  \right]
$. Calculation of $\mathbb{E}K^{2}\left(  \left\Vert X-x_{0}\right\Vert
/h\right)  $ is similar with $\mathbb{E}K\left(  \left\Vert X-x_{0}\right\Vert
/h\right)  $. In a multivariate setting, when $X$ is an $\mathbb{R}^{d}$
valued random variable, and the density of $X$ $f_{X}$ is smooth enough at
$x_{0}$ computations lead in many situations to :%
\begin{equation}
\mathbb{E}K\left(  \left\Vert X-x_{0}\right\Vert /h\right)  \sim c_{d}%
f_{X}\left(  x_{0}\right)  h^{d} \label{findim}%
\end{equation}
where $c_{d}$ denotes the volume of the unit ball in the space $\mathbb{R}%
^{d}$. The r.h.s. of the formula above may vary, depending on the support of
the distribution of $X$. However neither Lebesgue's measure or a counterpart
to $f_{X}$ may be defined when $X$ is valued in a Hilbert space for instance.
The classical notion of volume of a ball cannot be generalized to such spaces.
As a consequence when $X$ is a process, the density of $X$ at $x_{0}$ does not
make sense anymore. A major issue is then to compute the preceding expectation
without assuming that $f_{X}\left(  x_{0}\right)  $ exists. We consider the
following conditions on the kernel $K$ :

\begin{center}
$K$\textbf{ has compact support (say }$\left[  0,1\right]  $\textbf{), is
absolutely continuous and bounded above and below with }$K\left(  1\right)
>0$
\end{center}

These conditions hold for the naive kernel, $K\left(  u\right)  =1$ if and
only if $u\in\left[  0,1\right]  $. We do not seek minimal conditions on the
kernel here and the assumption above could certainly be alleviated but is
sufficient to carry out computations. Applying Fubini's Theorem is sufficient
to get rid of the density. Denoting $\mathbb{P}\left(  \left\Vert
X-x_{0}\right\Vert <h\right)  =F_{x_{0}}\left(  h\right)  $ we obtain :%
\begin{equation}
\mathbb{E}K\left(  \frac{\left\Vert X-x_{0}\right\Vert }{h}\right)  =F_{x_{0}%
}\left(  h\right)  \left[  K\left(  1\right)  -\int_{0}^{1}K^{\prime}\left(
s\right)  \frac{F_{x_{0}}\left(  hs\right)  }{F_{x_{0}}\left(  h\right)
}ds\right] \nonumber
\end{equation}
It is straightforward to see that the same method may yield the value of such
integrals as :%
\begin{equation}
\mathbb{E}\left[  \left\Vert X-x_{0}\right\Vert ^{p}K\left(  \frac{\left\Vert
X-x_{0}\right\Vert }{h}\right)  \right]  =h^{p}F_{x_{0}}\left(  h\right)
\left[  K\left(  1\right)  +\int_{0}^{1}\widetilde{K}_{p}\left(  s\right)
\frac{F_{x_{0}}\left(  hs\right)  }{F_{x_{0}}\left(  h\right)  }ds\right]
\label{exp-ker2}%
\end{equation}
with $\widetilde{K}_{p}\left(  s\right)  =-\left[  s^{p}K^{\prime}\left(
s\right)  +ps^{p-1}K\left(  s\right)  \right]  $ and the evaluation of the
expectation above essentially depends again on the small ball probability
$F_{x_{0}}\left(  \cdot\right)  $. When $X$ is a random function the behaviour
of $F_{x_{0}}$ at $0$ is crucial and determines the rate of convergence to
zero of the above expectation -what statisticians are truly interested in.

Assume that $F_{x_{0}}$ is regularly varying at zero with index $d$ (which is
usually true when $X$ is finite dimensional) then by definition $F_{x_{0}%
}\left(  h\right)  =Ch^{d}l\left(  h\right)  $ where $C$ is a constant, $l$ is
a slowly varying function at $0$ and $F_{x_{0}}\left(  hs\right)  /F_{x_{0}%
}\left(  h\right)  \rightarrow s^{d}$ when for $s>0$ and $h\rightarrow0$ which
yields $\mathbb{E}K\left(  \left\Vert X-x_{0}\right\Vert /h\right)  \sim
c_{d}h^{d}l\left(  h\right)  $ where $c_{d}$ depends only on $d$ and $K$.
Unfortunately when $X$ lies in a function space, the most classical examples
of $F\left(  h\right)  $ are not reguarly varying as will be seen below. But
however we notice for further purpose that the theory of regular variation is
of some help in this important special case.

Turning to $\mathbb{E}\left[  \left\langle X,e_{i}\right\rangle ^{2}K\left(
\left\Vert X\right\Vert /h\right)  \right]  $ which appears in the numerator
of (\ref{b}) (note that here $x_{0}$ does not appear anymore) is more tricky
and will not be done at this stage. This expectation is bounded above by
$\mathbb{E}\left[  \left\Vert X\right\Vert ^{2}K\left(  \left\Vert
X\right\Vert /h\right)  \right]  $ similar to (\ref{exp-ker2}) with $p=2$ but
this boudn is not sharp and no other equivalent could be derived from the
previous considerations.

\section{Main results}

We are ready to give the main results. This section is split in three parts.
In the first it is shown that the function $F_{x_{0}}\left(  \cdot\right)  $
which is crucial for evaluating the risk in model (\ref{reg-mod}) belongs to
the class $\Gamma_{0}$ of Gamma-varying functions in a quite general
framework. In the second we focus on the case of a Gaussian design. In the
third we use the properties of the class $\Gamma_{0}$ to derive upper and
lower bounds on the risk for (\ref{reg-mod}) and at a fixed point. A notable
fact is that the lower bound is degenerate : it is slower than any negative
power of $n.$ This may be seen as an ultimate symptom of the curse of
dimensionality. If $f$ and $g$ are two positive functions the notation
$f\preceq_{x}g$ means that $\lim_{u\rightarrow x}f\left(  u\right)  /g\left(
u\right)  \leq c$ for some positive constant $c.$

\subsection{Small ball probability of random functions are Gamma-varying}

This sub-section connects the two apparently distinct notions of probability
seen before : the class of small ball probabilities in $l_{2}$ and de Haan's
Gamma class of functions. Both families of functions are defined by their
local behaviour around $0.$ In what follows, the exponent $-1$ is strictly
reserved to denoting the generalized inverse of a function $f$ denoted
$f^{-1}$. Consequently in general $f^{-1}\neq1/f.$ Let us introduce the
function $\lambda\left(  \cdot\right)  $ which interpolates the $\lambda_{j}%
$'s a smooth way (which means that $\lambda\left(  j\right)  =\lambda_{j}$ for
all $j$ and $\lambda$ is $\mathbf{C}^{1}$).

Since our results rely on those of Lisfhits (1997) we recall now the
assumptions needed in this article. Let $G$ denote the (cumulative)
distribution function of $Z$ then we assume that there exists $b\in\left(
0,1\right)  ,$ $c_{1}>1$, $c_{2}\in\left(  0,1\right)  $ and $c_{3}>0$ such
that for $r<c_{3}$ :%
\begin{equation}
\mathbf{A}_{0}:\left\{
\begin{array}
[c]{c}%
G\left(  r\right)  \leq c_{1}G\left(  br\right) \\
G\left(  br\right)  \leq c_{2}G\left(  r\right) \\
\mathbb{E}Z^{3}<+\infty
\end{array}
\right.  \label{A0}%
\end{equation}

As mentioned in Lifhsits (1997) assumption $\mathbf{A}_{0}$ states that the
local behaviour at $0$ of $G$ is polynomial and $\mathbf{A}_{0}$ holds
whenever the density $g$ of $Z$ is regularly varying at $0$ with index
$\alpha>-1$. We also note that the assumption above holds for a large class of
classical positive distributions of $Z$ itself (Gamma, Beta...) or when
$Z=X^{2}$ with $X$ Gaussian, $X$ Laplace, Uniform or Student distributions for
instance. These considerations are of interest for the statistician in order
not to limit the approach to Gaussian models. Note that the assumption on the
convergence of the third order moment of $Z$ was alleviated in some recent
papers. We keep it here since it is general enough for our purpose.

When $\left(  Z_{i}\right)  _{i\in\mathbb{N}}$ is a sequence of random
variables whose cumulative distribution function $G$ is regularly varying at
$0$ with strictly positive index, the explicit form of the small ball
probability was derived for explicit sequences of log convex $\lambda\left(
\cdot\right)  $ by Dunker, Lifshits, Linde (1998). In particular they show
that when $\lambda_{i}=i^{-\beta}$ $\left(  \beta>1\right)  ,$ $\mathbb{P}%
\left(  \left\Vert X\right\Vert ^{2}<s\right)  \sim F_{1}\left(  s\right)  $
and that when $\lambda_{i}=\exp\left(  -i\right)  $ $\mathbb{P}\left(
\left\Vert X\right\Vert ^{2}<s\right)  \sim F_{2}\left(  s\right)  $ with :%
\begin{align}
F_{1}\left(  s\right)   &  =c_{1}s^{\left[  1+\beta\left(  2+c_{2}\right)
\right]  /\left(  2\beta-2\right)  }\exp\left(  -c_{3}s^{-1/\left(
\beta-1\right)  }\right) \label{smb1}\\
F_{2}\left(  s\right)   &  \sim c_{4}\left[  s^{1/3}\log\left(  1/s\right)
\right]  ^{-3/4}\exp\left(  -\left[  \log\left(  s\log1/s\right)  \right]
^{2}/4+\psi_{0}\left(  \log\left(  s\log1/s\right)  \right)  \right)
\label{smd2}%
\end{align}
where $\psi_{0}$ is a bounded function. Formula (\ref{smb1}) is proved as well
at page 269 in Lifshits (1995). Simple algebra proves that both functions on
the right hand side of (\ref{smb1}) and (\ref{smd2}) have all their
derivatives vanishing at $0.$ We notice that the r.h.s. of (\ref{smb1}) is
always flatter than the r.h.s. of (\ref{smd2}) which in turn will always be
flatter at $0$ than any polynomial function (like $c_{d}s^{d}$). However we
notice that the degree of flatness is directly connected with the rate of
decrease of the $\lambda_{i}$'s which quantifies, exactly like the l-numbers,
the accuracy of a finite-dimensional approximation of $X.$ We emphasize the
following Proposition, which will not be proved, on purpose.

\begin{proposition}
\label{oi}Both functions $F_{1}$ and $F_{2}$ defined above at (\ref{smb1}) and
(\ref{smd2}) belong to $\Gamma_{0}$ with respective auxiliary functions
$\rho_{1}\left(  s\right)  \sim s^{\beta/\left(  \beta-1\right)  }$ and
$\rho_{2}\left(  s\right)  \sim s\log\left(  1/s\right)  $ which both match
Proposition \ref{cabrel}.
\end{proposition}

The auxiliary functions $\rho_{1}$ and $\rho_{2}$ could be more precisely
computed but we only need equivalencies at this stage.\newline We are ready to
extend this fact to general sequences $\left(  \lambda_{i}\right)
_{i\in\mathbb{N}}$. Remind that the function $\gamma\left(  \cdot\right)  $
was defined implicitely at display (\ref{mu}). In words it is, up to sign, the
inverse of the first order derivative of the log-Laplace transform of
$S=\sum_{i=1}^{+\infty}\lambda_{i}Z_{i}$.

\begin{theorem}
\label{main} Let $S$ be defined by (\ref{sbpl2}) and set $\mathbb{P}\left(
S<s\right)  =F\left(  s\right)  $ the small ball probability of $S$ then
$F\in\Gamma_{0}$ with auxiliary function :%
\begin{equation}
\rho\left(  s\right)  =\frac{1}{\gamma\left(  s\right)  } \label{morcheeba}%
\end{equation}
and the representation (\ref{Dem}) may be rephrased only in terms of
$\gamma\left(  \cdot\right)  $ :%
\begin{equation}
\mathbb{P}\left(  S<r\right)  \underset{r\rightarrow0}{\sim}\frac{1}%
{\sqrt{2\pi}}\frac{\sqrt{-\gamma^{\prime}\left(  r\right)  }}{\gamma\left(
r\right)  }\exp\left[  -\int_{r}^{r_{0}}\gamma\left(  s\right)  ds\right]
\label{Lifs}%
\end{equation}
where $r_{0}=\mathbb{E}Z\cdot\sum_{j=1}^{+\infty}\lambda_{j}.$
\end{theorem}

Obviously the r.h.s. of (\ref{Lifs}) is mathematically the same object as the
r.h.s. of (\ref{Dem}). The "Gamma-varying version" of the r.h.s. is
$\sqrt{\rho^{\prime}\left(  r\right)  /\pi}\exp\left[  -\int_{r}^{r_{0}%
}ds/\rho\left(  s\right)  \right]  $. We believe however that this new version
is slightly more explicit and maybe more suited for statistical purposes. We
can take advantage as well of the properties of the class $\Gamma_{0}$ listed earlier.

The Theorem may be intuitively explained in view of Proposition \ref{skye}.
Indeed when $X$ lies in $\mathbb{R}^{d}$ and in a general context $F\left(
s\right)  \sim_{0}p_{d}\left(  s\right)  =c_{d}s^{d}$. The function $p_{d}$
has the following property : $p_{d}^{\left(  k\right)  }\left(  0\right)  =0$
whenever $k\neq d$. Consequently in an infinite dimensional space we can
expect that all the derivatives at $0$ should be null and this property is
recovered through Proposition \ref{skye}. A more geometric way to understand
this consists in considering the problem of the concentration of a probability
measure. Let $\mu$ be the measure associated with the random variable $X$.
Once again starting from $\mathbb{R}^{d}$ and letting $d$ increase -even if
this approach is not really fair- we see that $\mu$ must allocate a constant
mass of $1$ to a space whose dimension increases. Then $\mu$ gets more and
more diffuse, allowing fewer mass to balls and visiting rarely fixed points
such as $x_{0}$ (and their neighborhoods), resulting in a very flat small ball
probability function.

The following corollary provides some information about the rate of decrease
to zero of $F\left(  \cdot\right)  $ when an additional assumption is made on
$\rho$.

\begin{corollary}
\label{corrv}Assume that $\rho\left(  s\right)  =s^{\alpha}l\left(  s\right)
$ with $l\left(  \cdot\right)  $ slowly varying at $0$ (which just means that
$\rho$ is regularly varying at $0$ with index $\alpha\geq1$ and set :%
\begin{equation}
\mathbf{RV}_{+}:\alpha>1\mathrm{\ or\ }\mathbf{RV}_{1}:\rho\left(  s\right)
=sl\left(  s\right)  \mathrm{\ with\ }l\left(  s\right)  \succeq\log\left(
1/s\right)  . \label{rv}%
\end{equation}
If $\mathbf{RV}_{+}$ holds $\log\mathbb{P}\left(  S<r\right)  \preceq
c_{\alpha}r^{1-\alpha}$ and when $\mathbf{RV}_{1}$ holds $\log\mathbb{P}%
\left(  S<r\right)  \preceq-\varsigma\left(  r\right)  \log\left(  1/r\right)
$ for some $\varsigma\left(  r\right)  \rightarrow+\infty$ when $r\rightarrow
0$. In both preceding cases for all integer $p$ $\lim_{r\rightarrow0}%
r^{-p}\mathbb{P}\left(  S<r\right)  =0.$
\end{corollary}

This property fo the small ball probability has to be connected with property
(\ref{zaza}), is referred to as "rapid variation" at $0$ in the literature on
regular variations and may be compared or opposed with the regularly varying
situation discussed below (\ref{exp-ker2}). The assumptions $\mathbf{RV}_{+}$
and $\mathbf{RV}_{1}$ will be encountered again when addressing the case of
nonparametric regression. At last, note that for the auxiliary functions
$\rho_{1}$ and $\rho_{2}$ appearing at Proposition \ref{oi} and arising from
Dunker, Lifshits, Linde (1998) work we get $\rho_{1}\in\mathbf{RV}_{+}$ and
$\rho_{2}\in\mathbf{RV}_{1}$.

\textbf{Proof of Corollary \ref{corrv}:} We focus on the right hand side of
(\ref{Lifs}). First from $\sqrt{-\gamma^{\prime}\left(  r\right)  }%
/\gamma\left(  r\right)  =\sqrt{2\rho^{\prime}\left(  s\right)  }$ and the
properties of the auxiliary function $\rho$ at Proposition \ref{cabrel} we
have that $\rho^{\prime}\left(  s\right)  \rightarrow0$. Hence $\mathbb{P}%
\left(  S<r\right)  \leq\exp\left[  -\int_{r}^{r_{0}}1/\rho\left(  s\right)
ds\right]  $ for $r$ tending to $0$. From the direct part of Karamata's
theorem's (see Bingham, Goldie Teugels (1991), p.26) $\int_{r}^{r_{0}}%
1/\rho\left(  s\right)  ds\sim c_{\alpha}r^{1-\alpha}$ when $\alpha>1$ and
$\int_{r}^{r_{0}}1/\rho\left(  s\right)  ds=\varsigma\left(  r\right)
\log\left(  1/r\right)  $ with $\varsigma\left(  r\right)  \rightarrow+\infty$
when $r\rightarrow0$ (see display (1.5.8) in Bingham, Goldie Teugels (1991)).
Finally when $\mathbf{RV}_{+}$ or $\mathbf{RV}_{1}$ hold $-p\log r-\int
_{r}^{r_{0}}1/\rho\left(  s\right)  ds$ always tend to $-\infty$ whatever the
choice of $p$.

\begin{remark}
For the sake of completeness we point out the following fact which may be
misleading : indeed we started from $\mathbb{P}\left(  \left\Vert X\right\Vert
^{2}<r\right)  $ and the properties of this function may differ from those of
what may be intended as the true small ball probability $\mathbb{P}\left(
\left\Vert X\right\Vert ^{2}<r^{2}\right)  .$ It is not difficult to show that
if $F\in\Gamma_{0}$ with auxiliary function $\rho_{F}$ then $G\left(
r\right)  =F\left(  r^{2}\right)  $ belongs to $\Gamma_{0}$ as well with
auxiliary function $\rho_{G}$ defined by $\rho_{G}\left(  r\right)  =\rho
_{F}\left(  r^{2}\right)  /\left(  2r\right)  .$
\end{remark}

\subsection{Gaussian framework}

Assuming that $X$ is Gaussian, hence that $x_{i}$ are $\mathcal{N}\left(
0,1\right)  $ distributed provides a critical amount of extra information.
Indeed it is then possible to compute in a more explicit form :%
\begin{equation}
r=-\frac{\partial\log\Lambda\left(  \gamma\right)  }{\partial\gamma}=\sum
_{j}\frac{\lambda_{j}}{1+2\gamma\lambda_{j}} \label{log-lap-gauss}%
\end{equation}
which is the seminal equation linking $r$ and $\gamma$. We derive below an
explicit link between the $\lambda_{j}$'s and $\gamma\left(  \cdot\right)  $
or equivalently $\rho\left(  \cdot\right)  $. Under rather general assumptions
on the rate of decrease of the $\lambda_{j}$'s we obtain as well an upper
bound for the small ball probability which will be exploited in the next
subsection when investigating a lower bound for the regression.

\begin{proposition}
\label{eq}Assume that $X$ is Gaussian, that $\lambda\left(  \cdot\right)  $ is
a convex decreasing function and set $\varphi\left(  t\right)  =t\gamma\left(
t\right)  $. We have the following : There exists a fixed constant
$l\in\left[  \frac{2}{3},\frac{3}{2}\right]  $ such that for any
$\varepsilon>0$ and large enough $x$%
\[
\lambda\left(  x\left(  l+\varepsilon\right)  \right)  \leq\frac{1}%
{\gamma\left(  \varphi^{-1}\left(  x\right)  \right)  }=\rho\left(
\varphi^{-1}x\right)  \leq\lambda\left(  x\left(  l-\varepsilon\right)
\right)
\]
Besides when $\lambda\left(  x\right)  \succ_{\infty}\exp\left(  -x^{\alpha
}\right)  $ for some $\alpha>0,$ $F\left(  s\right)  \prec_{0}\exp\left[
-\left(  \log1/r\right)  ^{1+1/\alpha}\right]  .$ When $\lambda\left(
\cdot\right)  $ is explicitely known more precise relationships may be
derived. For instance when $\lambda_{a}\left(  x\right)  =cx^{-1-\nu}$ with
$c,\nu>0,$ $\gamma_{a}\left(  s\right)  \sim_{0}s^{-1-1/\nu}c^{1/\nu
}l^{1+1/\nu}$ and $l=\int_{0}^{+\infty}du/\left(  2+u^{1+\nu}\right)  $. When
$\lambda_{g}\left(  x\right)  =c\exp\left(  -\nu x\right)  ,$ with again
$c,\nu>0,$ $\gamma_{g}\left(  s\right)  \sim_{0}\log\left(  1/s\right)
/\left(  \nu s\right)  .$
\end{proposition}

\begin{remark}
The auxiliary fucntions $\rho_{a}=1/\gamma_{a}$ and $\rho_{g}=1/\gamma_{g}$
match respectively $\rho_{1}$ and $\rho_{2}$. Besides letting $\nu$ go to
infinity we see that, in a way $\rho_{g}$ may be viewed as a limit of
$\rho_{a}$. In fact $1/\log\left(  1/s\right)  $ echoes the degeneracy of
$1/\left(  \beta-1\right)  $. Proposition \ref{eq} modestly rediscover the
results of Dunker, Lifshits, Linde (1998). The upper bound of $F\left(
\cdot\right)  $ is close to the one obtained in Proposition \ref{corrv}. No
asumptions are needed on $\rho$ here but the distribution of $X$ is Gaussian.
\end{remark}

\textbf{Proof of Proposition \ref{eq} :}

We start from (\ref{log-lap-gauss}) and denote $a\left(  \cdot\right)
=1/\lambda\left(  \cdot\right)  $ which may be rewritten :%
\[
r=\sum_{j\geq1}\frac{1}{a\left(  j\right)  +2\gamma}%
\]

Let us set $J_{\gamma}=\inf\left\{  j:a\left(  j\right)  \geq\gamma\right\}  $
so that $a\left(  J_{\gamma}-1\right)  \leq\gamma\leq a\left(  J_{\gamma
}\right)  $%

\begin{align*}
\frac{J_{\gamma}}{3\gamma}+\sum_{j\geq J_{\gamma}+1}\frac{1}{a\left(
j\right)  +2\gamma}  &  \leq\sum_{j\geq1}\frac{1}{a\left(  j\right)  +2\gamma
}\leq\frac{J_{\gamma}}{2\gamma}+\sum_{j\geq J_{\gamma}+1}\frac{1}{a\left(
j\right)  +2\gamma}\\
\frac{J_{\gamma}}{3\gamma}+\frac{1}{3}\sum_{j\geq J_{\gamma}}\frac{1}{a\left(
j\right)  }  &  \leq r\leq\frac{J_{\gamma}}{2\gamma}+\sum_{j\geq J_{\gamma}%
}\frac{1}{a\left(  j\right)  }\\
\frac{2}{3}  &  \leq\frac{r\gamma}{J_{\gamma}}\leq\left(  0.5+\left(
J_{\gamma}^{-1}+1\right)  \right)
\end{align*}
The convexity of $\lambda,$ hence of $a$ yields $\sum_{j\geq J}1/a\left(
j\right)  \leq\left(  J_{\gamma}+1\right)  /a\left(  J_{\gamma}\right)  $ (see
Cardot, Mas, Sarda (2007)) hence for $r\downarrow0$ that is for large
$\gamma\left(  r\right)  $ :%
\begin{equation}
\frac{2}{3}\leq\frac{r\gamma}{J_{\gamma}}\leq\frac{3}{2} \label{sc}%
\end{equation}
Now consider the function (of the variable $\gamma$) : $d\left(
\gamma\right)  =\gamma r\left(  \gamma\right)  /a^{-1}\left(  \gamma\right)  $
is decreasing at least when $\gamma\uparrow+\infty$ since :%
\[
d^{\prime}\left(  \gamma\right)  =\frac{\left(  r\left(  \gamma\right)
-\gamma\sigma^{2}\right)  }{a^{-1}\left(  \gamma\right)  }-\frac{\gamma
r\left(  \gamma\right)  }{\left[  a^{-1}\left(  \gamma\right)  \right]
^{2}\cdot a^{\prime}\left(  a^{-1}\left(  \gamma\right)  \right)  }%
\]
is negative for $r\downarrow0$ that is for large $\gamma$. Finally we get
$r\gamma/a^{-1}\left(  \gamma\right)  \rightarrow l\in\left[  \frac{2}%
{3},\frac{3}{2}\right]  $ when $\gamma$ tends to $+\infty$. The first
statement of the Theorem is a consequence of the latter limit. As a
consequence when $\lambda\left(  x\right)  \succ_{\infty}\exp\left(
-x^{\alpha}\right)  $, $a\left(  x\right)  \prec_{\infty}\exp\left(
x^{\alpha}\right)  $ and finally $r\gamma\succ_{0}\left(  \log\gamma\right)
^{1/\alpha}\geq\left(  \log1/r\right)  ^{1/\alpha}$. At last we get
$\gamma\left(  r\right)  \succ_{0}\frac{1}{r}\left(  \log1/r\right)
^{1/\alpha}$ which finally yields for some constant $c$%
\begin{align*}
\exp\left[  -\int_{r}^{r_{0}}\gamma\left(  s\right)  ds\right]   &  \leq
c\exp\left[  \int_{r}^{r_{0}}\left(  \log1/s\right)  ^{1/\alpha}d\left(
\log1/s\right)  \right] \\
&  =c\exp\left[  -\left(  \log1/r\right)  ^{1+1/\alpha}\right]
\end{align*}

The last sentence of the Theorem , when the function $\lambda$ is known, is
easily derived by noting that $\sum_{j\geq1}1/\left(  a\left(  j\right)
+2\gamma\right)  \sim_{\gamma\rightarrow+\infty}\int_{0}^{+\infty}\frac
{dx}{a\left(  x\right)  +2\gamma}.$

\subsection{Upper and lower bound in regression for functional data}

We fix once and for all the assumptions considered in what follows. These
assumptions appear in addition to those considered in the previous sections.
Remind that if $g$ is some function defined on $\mathcal{H}$ and with values
in $\mathbb{R}$ the first order Fr\'{e}chet-derivative of $g$ at $x_{0}$ (its
infinite-dimensional gradient) may be identified with an element of
$\mathcal{H}$. The second order derivative $g^{\prime\prime}\left(
x_{0}\right)  $ is identified with a symmetric operator from $\mathcal{H}$ to
$\mathcal{H}.$

\textbf{Assumptions on the distribution of }$X.$ The random element $X$ is
centered and in the development (\ref{KL}) the $x_{i}$'s are independent. We
have $P_{X-x_{0}}\ll P_{X}$ with $f_{x_{0}}=dP_{X-x_{0}}/dP_{X}$ such that
$f_{x_{0}}\left(  0\right)  >0$, $f^{\prime}\left(  x_{0}\right)
\in\mathcal{H}$ exists and the second order derivative of $f_{x_{0}}$ denoted
$f_{x_{0}}^{\prime\prime}\left(  x\right)  $ is for all $x$ in a neighborhood
of $x_{0}$ a bounded linear operator from $\mathcal{H}$ to $\mathcal{H}$.
Denote $\partial_{i}f_{x_{0}}=\left\langle f^{\prime}\left(  x_{0}\right)
,e_{i}\right\rangle $ where $e_{i}$ is one of the eigenvectors appearing in
(\ref{KL}). Besides we assume that for all $i$ the density of the margins
$\left\langle X,e_{i}\right\rangle $ is symmetric.

\textbf{Assumptions on the regression function}. Assume that $r$ has first and
second order derivative at $x_{0}$. We denote $\partial_{i}r_{x_{0}%
}=\left\langle r^{\prime}\left(  x_{0}\right)  ,e_{i}\right\rangle $ and
$\partial_{ii}^{2}r_{x_{0}}=\left\langle r^{\prime}\left(  x_{0}\right)
\left(  e_{i}\right)  ,e_{i}\right\rangle $ and assume as well that :%
\[
\sum_{i=1}^{+\infty}\lambda_{i}\partial_{i}r_{x_{0}}\partial_{i}f_{x_{0}}\neq0
\]

At this point a discussion on the assumptions related to the distribution of
$X$ is needed. Take for instance the case of a gaussian $X.$ Chapter 9 and 10
in Lifshits (1995) are clear about these issues (see more specifically
p.102-107.) It is possible to shift the assumptions on the regularity of $X$
to conditions on the regularity of $x_{0}$. First in order to define
$f^{\prime}\left(  x_{0}\right)  $ we need to assume that $x_{0}=\left(
m_{i}\right)  _{1\leq i}$ belongs to the kernel of $X$ that is $\sum m_{i}%
^{2}/\lambda_{i}<+\infty.$ Then for any $u=\left(  u_{i}\right)  _{1\leq i}$
in $\mathcal{H}$ $f_{x_{0}}\left(  u\right)  =\exp\left(  -\sum_{i\geq1}%
m_{i}^{2}/2\lambda_{i}+\sum_{i\geq1}u_{i}m_{i}/2\lambda_{i}\right)  $and
$\partial_{i}f_{x_{0}}=\partial_{i}f_{x_{0}}\left(  u\right)  =\left(
m_{i}/2\lambda_{i}\right)  f_{x_{0}}\left(  u\right)  $. We are interested in
th smoothness of these functions at $0.$ From $\left\vert f_{x_{0}}\left(
u\right)  -f_{x_{0}}\left(  0\right)  \right\vert \preceq\left\Vert
u\right\Vert \left(  \sum_{i\geq1}m_{i}^{2}/\lambda_{i}^{2}\right)  ^{1/2}$.
And the finiteness of the latter series is subject to a condition of decay on
the coefficients $m_{i}$'s.

It turns out that the gaussian framework may be generalized to some other
distributions. It suffices to consider a symmetric random variable $U$ such
that $\mathbb{E}U=0$ and $\mathbb{V}U=1$ with smooth density function at
$0.$Then introduce the scaled $Z_{i}=\lambda_{i}U^{2}$ in order to derive a
new development (\ref{KL}). Examples for $U$ are the uniform distribution on
$\left[  -a,a\right]  $, $a>0$, shifted Beta distributions...

\subsubsection{Upper bound}

The reader was left with Lemma \ref{risk}. In view of the results of this
section we are in a position to simplify some computations. Turning to the
local moments defined at display (\ref{exp-ker2}), from properties of
functions in $\Gamma_{0}$ and specifically (\ref{zaza}) we get%
\begin{equation}
\mathbb{E}\left[  \left\Vert X-x_{0}\right\Vert ^{p}\right]  K\left(
\frac{\left\Vert X-x_{0}\right\Vert }{h}\right)  \sim K\left(  1\right)
h^{p}F_{x_{0}}\left(  h\right)  \sim K\left(  1\right)  h^{p}f_{x_{0}}\left(
0\right)  8F\left(  h\right)  \label{gvexpmom}%
\end{equation}
We see again that the representation theorem of the preceding section is of
some help to simplify our calculations. We mention for immediate purpose that
the derivation of both fromula above leads as well to :%
\begin{equation}
\mathbb{E}K^{2}\left(  \frac{\left\Vert X-x_{0}\right\Vert }{h}\right)  \sim
K^{2}\left(  1\right)  F_{x_{0}}\left(  h\right)  \label{sq-kern}%
\end{equation}
Let the local moments of order $1$ and $2$ of $X$ at $x_{0}$ be respectively
defined by~:%
\begin{align}
\mathcal{M}_{K,1}\left(  x_{0}\right)   &  =\mathbb{E}\left[  \left(
X-x_{0}\right)  K\left(  \frac{\left\Vert X-x_{0}\right\Vert }{h}\right)
\right] \label{noel}\\
\mathcal{M}_{K,2}\left(  x_{0}\right)   &  =\mathbb{E}\left[  \left(
X-x_{0}\right)  \otimes\left(  X-x_{0}\right)  \right]  K\left(
\frac{\left\Vert X-x_{0}\right\Vert }{h}\right)  . \label{smet}%
\end{align}
Formula (\ref{smet}) may be explicited. First let $u$ and $v$ be two points in
the vector space then $u\otimes v$ is a linear operator defined by $\left[
u\otimes v\right]  \left(  x\right)  =\left\langle v,x\right\rangle y$. Now
erasing $x_{0}$ and $K\left(  \left\Vert X-x_{0}\right\Vert /h\right)  $ gives
the usual covariance operator of $X$ for a centered $X.$ The special
covariance operator $\mathcal{M}_{K,2}\left(  x_{0}\right)  $ is obtained by
shifting and smoothing $X$ around $x_{0}$. Note that $\mathcal{M}_{K,1}$
belongs to $\mathcal{H}$ and $\mathcal{M}_{K,2}$ is a linear trace-class
operator acting from and onto $\mathcal{H}$. We refer to M\"{u}ller and Yan
(2001) for some statistical results on local moments for finite-dimensional
random variables and to Mas (2008) for some related results dealing with
(\ref{smet}) and where random functions and small ball problems appear.

The next Proposition completes Lemma \ref{risk}.

\begin{proposition}
\label{risk2}For the variance part of the risk the equivalence holds
$\mathcal{V}_{n}\left(  x_{0}\right)  \sim\sigma_{\varepsilon}^{2}/\left[
f_{x_{0}}\left(  0\right)  nF\left(  h\right)  \right]  $. For the bias part
we just provide an approximate rate :%
\[
\mathbf{c}^{-}\rho^{6}\left(  h\right)  \leq\mathcal{B}_{n}\left(
x_{0}\right)  \leq\mathbf{c}^{+}h^{4}%
\]
where $\mathbf{c}^{+}$ and $\mathbf{c}^{-}$ depend only on $f_{x_{0}}\left(
0\right)  ,$ $f^{\prime}\left(  x_{0}\right)  ,$ $r^{\prime}\left(
x_{0}\right)  $ and $r^{\prime\prime}\left(  x_{0}\right)  $. When$\rho\left(
h\right)  \succeq h^{m}$ for some $m$ $\mathcal{B}_{n}\left(  x_{0}\right)  $
decreases to $0$ at most and at least at a polynomial rate.
\end{proposition}

The problem here is to ensure a rough control of $\mathcal{B}_{n}\left(
x_{0}\right)  $. As will be seen soon $\rho^{6}\left(  h\right)  $ turns out
to be regularly varying in most cases and decays to zero at a polynomial rate.
The unusual framework (namely with distributions in the class $\Gamma_{0}$)
motivates to prove the reader that $\mathcal{B}_{n}\left(  x_{0}\right)  $
does not reach an unusual rate of decreaes to $0$ (namely exponential). And
the lower bound we obtain for this specific estimator justifies the conditions
under which the minimax lower bound is going to be derived.

\textbf{Proof of Proposition (\ref{risk2}) :}

We start with $\mathcal{V}_{n}\left(  x_{0}\right)  $. It is simple to see
that $\mathcal{V}_{n}\left(  x_{0}\right)  =n\sigma_{\varepsilon}%
^{2}\mathbb{E}\omega_{1,n}^{2}$with :%
\[
\omega_{1,n}=\frac{K\left(  \left\Vert X_{1}-x_{0}\right\Vert /h\right)
}{\sum_{i=1}^{n}K\left(  \left\Vert X_{i}-x_{0}\right\Vert /h\right)  }%
\]
Computations like those carried out in Ferraty, Mas, Vieu (2007) show that :%
\[
\mathbb{E}\omega_{1,n}^{2}\sim\mathbb{E}K^{2}\left(  \left\Vert X-x_{0}%
\right\Vert /h\right)  /\left[  n\mathbb{E}K\left(  \left\Vert X-x_{0}%
\right\Vert /h\right)  \right]  ^{2}%
\]
hence that (see (\ref{gvexpmom})) $\mathcal{V}_{n}\left(  x_{0}\right)
\sim\frac{\sigma_{\varepsilon}^{2}}{n}\frac{\mathbb{E}K^{2}\left(  \left\Vert
X-x_{0}\right\Vert /h\right)  }{\left[  \mathbb{E}K\left(  \left\Vert
X-x_{0}\right\Vert /h\right)  \right]  ^{2}}\sim\frac{\sigma_{\varepsilon}%
^{2}}{nF_{x_{0}}\left(  h\right)  }$which yields the desired result by
(\ref{eqshsbp}).

We should now deal with $\mathcal{B}_{n}\left(  x_{0}\right)  $ defined at
(\ref{bias}). Ferraty, Mas, Vieu (2007 show that $\mathcal{B}_{n}\left(
x_{0}\right)  $ is well approximated by :%
\begin{equation}
\left[  \frac{\mathbb{E}\left(  r\left(  X\right)  -r\left(  x_{0}\right)
\right)  K\left(  \left\Vert X-x_{0}\right\Vert /h\right)  }{\mathbb{E}%
K\left(  \left\Vert X_{1}-x_{0}\right\Vert /h\right)  }\right]  ^{2}
\label{b2}%
\end{equation}
A Taylor development at order 2 shows that :%
\[
\mathbb{E}\left(  r\left(  X\right)  -r\left(  x_{0}\right)  \right)  K\left(
\left\Vert X-x_{0}\right\Vert /h\right)  =\left\langle r^{\prime}\left(
x_{0}\right)  ,\mathcal{M}_{K,1}\left(  x_{0}\right)  \right\rangle
+\mathrm{tr}\left[  r^{\prime\prime}\left(  \xi\right)  \mathcal{M}%
_{K,2}\left(  x_{0}\right)  \right]  /2
\]
where $r^{\prime}$ and $r^{\prime\prime}$ stands for the first and second
order G\^{a}teaux derivative of $r$ at $x_{0}$ and $\xi=\theta X+\left(
1-\theta\right)  x_{0}$ for some random $\theta\in\left(  0,1\right)  .$ We
first deal with the first order term $\left\langle r^{\prime}\left(
x_{0}\right)  ,\mathcal{M}_{K,1}\left(  x_{0}\right)  \right\rangle $. From
$\mathbb{P}_{X-x_{0}}\ll\mathbb{P}_{X}$ we see that :%
\begin{align*}
\mathbb{E}\left[  \left(  X-x_{0}\right)  K\left(  \frac{\left\Vert
X-x_{0}\right\Vert }{h}\right)  \right]   &  =\mathbb{E}\left[  Xf_{x_{0}%
}\left(  X\right)  K\left(  \frac{\left\Vert X\right\Vert }{h}\right)  \right]
\\
&  =\mathbb{E}\left[  X\left[  f_{x_{0}}\left(  X\right)  -f_{x_{0}}\left(
0\right)  \right]  K\left(  \frac{\left\Vert X\right\Vert }{h}\right)
\right]  +f_{x_{0}}\left(  0\right)  \mathbb{E}\left[  XK\left(
\frac{\left\Vert X\right\Vert }{h}\right)  \right]
\end{align*}
We assumed that $f_{i}$ the density of $x_{i}$ is symmetric. This yields for
all $i$ :%
\[
\mathbb{E}\left[  \left\langle X,e_{i}\right\rangle K\left(  \frac{\left\Vert
X\right\Vert }{h}\right)  \right]  =\mathbb{E}\left[  x_{i}K\left(
\frac{\left\Vert X\right\Vert }{h}\right)  \right]  =\int\left(  \int
tK\left(  \sqrt{t^{2}+s^{2}}/h\right)  f_{i}\left(  t\right)  dt\right)
g_{i}\left(  s\right)  ds=0
\]
where $g_{i}$ si the density of $\sum_{j=1,j\neq i}^{+\infty}\lambda_{j}%
x_{j}^{2}$. Now%
\[
\mathbb{E}\left[  X\left[  f_{x_{0}}\left(  X\right)  -f_{x_{0}}\left(
0\right)  \right]  K\left(  \frac{\left\Vert X\right\Vert }{h}\right)
\right]  =\mathbb{E}\left[  X\left\langle f^{\prime}\left(  x_{0}\right)
,X\right\rangle K\left(  \frac{\left\Vert X\right\Vert }{h}\right)  \right]
+R_{n}%
\]
where $R_{n}$ involves the second order derivative of $f_{x_{0}}$ and will be
neglected. Then denoting $\partial_{i}r_{x_{0}}=\left\langle r^{\prime}\left(
x_{0}\right)  ,e_{i}\right\rangle $%
\begin{align*}
\left\langle r^{\prime}\left(  x_{0}\right)  ,\mathcal{M}_{K,1}\left(
x_{0}\right)  \right\rangle  &  \sim\sum_{i=1}^{+\infty}\partial_{i}r_{x_{0}%
}\mathbb{E}\left[  \left\langle X,e_{i}\right\rangle \left\langle f^{\prime
}\left(  x_{0}\right)  ,X\right\rangle K\left(  \frac{\left\Vert X\right\Vert
}{h}\right)  \right] \\
&  =\sum_{i=1}^{+\infty}\partial_{i}r_{x_{0}}\mathbb{E}\left[  \left\langle
X,e_{i}\right\rangle \left(  \sum_{j=1}^{+\infty}\partial_{j}f_{x_{0}%
}\left\langle X,e_{j}\right\rangle \right)  K\left(  \frac{\left\Vert
X\right\Vert }{h}\right)  \right]
\end{align*}
Arguments based on the symmetry of the density of the $\left\langle
X,e_{i}\right\rangle $ lead to cancelling $\mathbb{E}\left[  \left\langle
X,e_{i}\right\rangle \left\langle X,e_{j}\right\rangle K\left(  \frac
{\left\Vert X\right\Vert }{h}\right)  \right]  $ for $i\neq j$ and :%
\[
\left\langle r^{\prime}\left(  x_{0}\right)  ,\mathcal{M}_{K,1}\left(
x_{0}\right)  \right\rangle \sim\sum_{i=1}^{+\infty}\partial_{i}r_{x_{0}%
}\partial_{i}f_{x_{0}}\mathbb{E}\left[  \left\langle X,e_{i}\right\rangle
^{2}K\left(  \frac{\left\Vert X\right\Vert }{h}\right)  \right]
\]
Similar calculations show that :%
\[
\mathrm{tr}\left[  r^{\prime\prime}\left(  x_{0}\right)  \mathcal{M}%
_{K,2}\left(  x_{0}\right)  \right]  \sim f_{x_{0}}\left(  0\right)
\sum_{i=1}^{+\infty}\partial_{ii}^{2}r_{x_{0}}\mathbb{E}\left[  \left\langle
X,e_{i}\right\rangle ^{2}K\left(  \frac{\left\Vert X\right\Vert }{h}\right)
\right]
\]
where $\partial_{ii}^{2}r_{x_{0}}=\left\langle r^{\prime\prime}\left(
x_{0}\right)  e_{i},e_{i}\right\rangle $ and finally denoting $\mathbb{E}%
\left[  \left\langle X,e_{i}\right\rangle ^{2}K\left(  \left\Vert X\right\Vert
/h\right)  \right]  =\mathbf{v}_{i}\left(  h\right)  $%
\[
\mathbb{E}\left(  r\left(  X\right)  -r\left(  x_{0}\right)  \right)  K\left(
\left\Vert X-x_{0}\right\Vert /h\right)  \sim\sum_{i=1}^{+\infty}\left(
\frac{\partial_{ii}^{2}r_{x_{0}}f_{x_{0}}\left(  0\right)  }{2}+\partial
_{i}r_{x_{0}}\partial_{i}f_{x_{0}}\right)  \mathbf{v}_{i}\left(  h\right)
\]

We can confine now derive an upper bound. Indeed for $h\downarrow0$,%
\begin{align*}
\left\vert \mathbb{E}\left(  r\left(  X\right)  -r\left(  x_{0}\right)
\right)  K\left(  \left\Vert X-x_{0}\right\Vert /h\right)  \right\vert  &
\leq2\sup_{i}\left\{  \frac{\partial_{ii}^{2}r_{x_{0}}f_{x_{0}}\left(
0\right)  }{2}+\left\vert \partial_{i}r_{x_{0}}\partial_{i}f_{x_{0}%
}\right\vert \right\}  \sum_{i=1}^{+\infty}\mathbf{v}_{i}\left(  h\right) \\
&  \leq2\sup_{i}\left\{  \frac{\partial_{ii}^{2}r_{x_{0}}f_{x_{0}}\left(
0\right)  }{2}+\left\vert \partial_{i}r_{x_{0}}\partial_{i}f_{x_{0}%
}\right\vert \right\}  \mathbb{E}\left[  \left\Vert X\right\Vert ^{2}K\left(
\left\Vert X\right\Vert /h\right)  \right] \\
&  \leq2K\left(  1\right)  \sup_{i}\left\{  \frac{\partial_{ii}^{2}r_{x_{0}%
}f_{x_{0}}\left(  0\right)  }{2}+\left\vert \partial_{i}r_{x_{0}}\partial
_{i}f_{x_{0}}\right\vert \right\}  h^{2}F\left(  h\right)
\end{align*}
This together with (\ref{b2}) and (\ref{gvexpmom}) leads to the upper bound of
the Proposition with $\mathbf{c}^{+}=8\left[  \sup_{i}\left\{  \partial
_{ii}^{2}r_{x_{0}}f_{x_{0}}\left(  0\right)  /2+\left\vert \partial
_{i}r_{x_{0}}\partial_{i}f_{x_{0}}\right\vert \right\}  /f_{x_{0}}\left(
0\right)  \right]  ^{2}$.

We turn to the lower bound. Since $r^{\prime\prime}\left(  \xi\right)
\mathcal{M}_{K,2}\left(  x_{0}\right)  $ is a positive operator we may confine
ourselves to the first term. For simplicity we will make calculations with the
naive kernel for $K$ and with a modified norm. In fact we will take
$\left\Vert X\right\Vert =\left\vert \left\langle X,e_{i}\right\rangle
\right\vert +\left\vert \sum_{j\neq i}\left\langle X,e_{j}\right\rangle
\right\vert =\left\vert \left\langle X,e_{i}\right\rangle \right\vert +Z_{i}$.
Let $0<\underline{c}<\overline{c}$ be two constants :%
\begin{align*}
\mathbb{E}\left[  \left\langle X,e_{i}\right\rangle ^{2}K\left(  \left\Vert
X\right\Vert /h\right)  \right]   &  \geq\lambda_{i}\underline{c}\rho
^{2}\left(  h\right)  \mathbb{P}\left(  \sqrt{\lambda_{i}\underline{c}}%
\rho\left(  h\right)  \leq\left\vert \left\langle X,e_{i}\right\rangle
\right\vert \leq\sqrt{\lambda_{i}\overline{c}}\rho\left(  h\right)
,\left\vert \left\langle X,e_{i}\right\rangle \right\vert +Z_{i}\leq h\right)
\\
&  \geq\lambda_{i}\underline{c}\rho^{2}\left(  h\right)  \mathbb{P}\left(
\left\vert \left\langle X,e_{i}\right\rangle \right\vert \in\sqrt{\lambda_{i}%
}\rho\left(  h\right)  \left[  \sqrt{\underline{c}},\sqrt{\overline{c}%
}\right]  \right)  \mathbb{P}\left(  Z_{i}\leq h-\sqrt{\overline{c}\lambda
_{i}}\rho\left(  h\right)  \right) \\
&  \geq\lambda_{i}\underline{c}\rho^{2}\left(  h\right)  \mathbb{P}\left(
\left\vert \left\langle X,e_{i}\right\rangle \right\vert \in\sqrt{\lambda_{i}%
}\rho\left(  h\right)  \left[  \sqrt{\underline{c}},\sqrt{\overline{c}%
}\right]  \right)  \mathbb{P}\left(  \left\Vert X\right\Vert \leq
h-\sqrt{\overline{c}\lambda_{i}}\rho\left(  h\right)  \right)
\end{align*}
where the probabilities were split because $\left\vert \left\langle
X,e_{i}\right\rangle \right\vert $ and $Z_{i}$ are independent. Consider first
$\mathbb{P}\left(  \left\vert \left\langle X,e_{i}\right\rangle \right\vert
/\sqrt{\lambda_{i}}\in\rho\left(  h\right)  \left[  \sqrt{\underline{c}}%
,\sqrt{\overline{c}}\right]  \right)  \geq\mathbf{c}\rho\left(  h\right)  $
where $\mathbf{c}$ is some constant independent of $i$ if the distribution of
all $\left\vert \left\langle X,e_{i}\right\rangle \right\vert /\sqrt
{\lambda_{i}}$ is bounded below in a neighborhood of $0$ which will be assumed
here (it is true when $X$ is gaussian). Then lastly :%
\[
F\left(  h-\sqrt{\overline{c}\lambda_{i}}\rho\left(  h\right)  \right)  \geq
F\left(  h-\sqrt{\overline{c}\lambda_{1}}\rho\left(  h\right)  \right)
\]
which yields :%
\[
\mathbb{E}\left[  \left\langle X,e_{i}\right\rangle ^{2}K\left(  \left\Vert
X\right\Vert /h\right)  \right]  \geq\lambda_{i}\underline{c}\mathbf{c}%
\rho^{3}\left(  h\right)  \mathbb{P}\left(  \left\Vert X\right\Vert \leq
h-\sqrt{\overline{c}\lambda_{1}}\rho\left(  h\right)  \right)
\]
From $F\left(  h-\sqrt{\overline{c}\lambda_{1}}\rho\left(  h\right)  \right)
/F\left(  h\right)  \rightarrow\exp\left(  -\sqrt{\overline{c}\lambda_{1}%
}\right)  $ we get that for a $h$ close enough to $0$ $\mathbb{E}\left[
\left\langle X,e_{i}\right\rangle ^{2}K\left(  \left\Vert X\right\Vert
/h\right)  \right]  \geq\lambda_{i}\mathbf{c}^{\prime\prime}\rho^{3}\left(
h\right)  F\left(  h\right)  $ where $\mathbf{c}^{\prime\prime}$ does not
depend on $n$ or on $i$. Finally we get :%
\begin{align*}
&  \left[  \frac{\mathbb{E}\left(  r\left(  X\right)  -r\left(  x_{0}\right)
\right)  K\left(  \left\Vert X-x_{0}\right\Vert /h\right)  }{\mathbb{E}%
K\left(  \left\Vert X_{1}-x_{0}\right\Vert /h\right)  }\right]  ^{2}\\
&  \geq\rho^{6}\left(  h\right)  \left[  \mathbf{c}^{\prime\prime}\sum
_{i=1}^{+\infty}\lambda_{i}\partial_{i}r_{x_{0}}\partial_{i}f_{x_{0}}%
/f_{x_{0}}\left(  0\right)  \right]  ^{2}\\
&  =\mathbf{c}^{-}\rho^{6}\left(  h\right)
\end{align*}
since $\sum_{i=1}^{+\infty}\lambda_{i}\partial_{i}r_{x_{0}}\partial
_{i}f_{x_{0}}\neq0.$

\subsubsection{Lower Bound}

From the preceding subsection the optimal risk for the kernel estimate is
obtained by selecting an $h$ balancing the trade-off between variance and
bias. Imagine that we had found in Proposition \ref{risk2} a result such as
$\mathcal{B}_{n}\left(  x_{0}\right)  \asymp F^{\kappa}\left(  h\right)  $ for
some $\kappa>0$. Then the optimal bandwidth would stem from $n^{-1}\asymp
F^{1+\kappa}\left(  h\right)  $ leading to a $\mathcal{R}_{n}\asymp
n^{-\kappa/\left(  1+\kappa\right)  }$ which would contradict the initial
claim of degenerate rate for the risk. This explains why we spend some energy
in delivering the lower bound on $\mathcal{B}_{n}\left(  x_{0}\right)  $ in
Proposition \ref{risk2}. As will be seen now when $r$ belongs to a class large
enough to inherit classical approximation features, $\mathcal{R}_{n}$ cannot
decrease at a ploynomial rate. What we mean by classical approximation
features is explicited now.

Let $\mathcal{E}_{p}$ denote any class of $\mathbb{R}$-valued functions
defined on $\mathcal{H}$ such that :%
\begin{equation}
\sup_{r\in\mathcal{E}_{p}}\mathcal{B}_{n}\left(  x_{0}\right)  \preceq h^{2p}
\label{reg-class}%
\end{equation}
For instance $\mathcal{E}_{p}$ may be the class of H\"{o}lder functions of
order $p\in\left]  0,1\right[  $. When $\mathcal{E}_{p}$ is the class of
function which have two derivatives at $x_{0}$ we see from Proposition
\ref{risk2} that (\ref{reg-class}) holds for some $p>0$ when $\rho\left(
h\right)  \geq h^{m}$ $m>1$. Optimizing the bias-variance trade-off in the
risk leads to choosing an $h$ such that $\sup_{r\in\mathcal{E}_{p}}%
\mathcal{B}_{n}\left(  x_{0}\right)  =\mathcal{V}_{n}\left(  x_{0}\right)  .$
The next Lemma deals with this issue.

\begin{lemma}
Assume that $X$ is gaussian and $\lambda\left(  x\right)  \succ_{\infty}%
\exp\left(  -x^{\alpha}\right)  $ for some $\alpha>0$. Let $c^{\ast}$ be some
constant and $h^{\ast}$ be the solution of the functional equation :%
\begin{equation}
\frac{1}{n}=c^{\ast}h^{2p}F\left(  h\right)  \label{opt-band}%
\end{equation}
then $n^{\beta}/\left(  nF\left(  h^{\ast}\right)  \right)  \rightarrow
+\infty$ for any $\beta>0.$When $X$ is non Gaussian but satisfies the
assumptions (\ref{A0}) and (\ref{rv}) the same conclusion holds.
\end{lemma}

\textbf{Proof of the Lemma :} Only the case $0<\beta<1$ has to be
investigated. When $X$ is Gaussian the lemma is easily derived from
Proposition \ref{eq} since it was proved that $F\left(  h\right)  \prec
_{0}\exp\left[  -\left(  \log1/h\right)  ^{1+1/\alpha}\right]  $ holds. When
$X$ is not gaussian and $\mathbf{RV}_{1}$ holds the proof of Corollary
\ref{corrv} shows that $\beta\log n+2p\log h^{\ast}>\beta\varsigma\left(
h^{\ast}\right)  \log\left(  1/h^{\ast}\right)  -2p\log h^{\ast}$where
$\varsigma\left(  h^{\ast}\right)  $ tends to $+\infty$ when $h^{\ast}$ tends
to $0.$ When $\mathbf{RV}_{+}$ holds the proof is the same with $c_{\alpha
}\left(  h^{\ast}\right)  ^{1-\alpha}$ instead of $\varsigma\left(  h^{\ast
}\right)  \log\left(  1/h^{\ast}\right)  $.\medskip

Now our approach to derive lower bounds for the minimax risk follows
Tsybakov's scheme (see Tsybakov (2004)) : we construct two models $r_{0}$ and
$r_{1}$ far enough from each other but such that the Hellinger distance
between the two models is bounded. Let $p_{\varepsilon}$ stand for the density
of $\varepsilon.$ Assume that for some constant $p_{\ast}$ and for all
$y\in\mathbb{R}$,%
\begin{equation}%
{\displaystyle\int_{\mathbb{R}}}
\left[  \sqrt{p_{\varepsilon}\left(  t\right)  }-\sqrt{\left[  p_{\varepsilon
}\left(  t+y\right)  \right]  }\right]  ^{2}dt\leq p_{\ast}y^{2}
\label{hyp-dens-eps}%
\end{equation}

This assumption is general and appears in Tsybakov's book. It holds under
smoothness assumptions on $p_{\varepsilon}$. We comment it briefly. If
$\Lambda\left(  y\right)  $ denotes the left hand side in the display above
$\Lambda\left(  y\right)  \leq2$ for all $y$ and we just need to study
$\Lambda$ on a compact neighborhood around $0$ (up to a rescaling through the
constant $p_{\ast}$). We see that $\Lambda\left(  0\right)  =0$ and
$\Lambda^{\prime}\left(  y\right)  =-\int p_{\varepsilon}^{\prime}\left(
t+y\right)  \sqrt{p_{\varepsilon}\left(  t\right)  /p_{\varepsilon}\left(
t+y\right)  }dt$ whenever $p_{\varepsilon}$ is smooth enough hence
$\Lambda^{\prime}\left(  0\right)  =0$. Under accurate conditions on
$\Lambda^{\prime\prime}$, $\Lambda\left(  y\right)  \leq p_{\ast}y^{2}$ will
hold around $0$ hence everywhere.

\begin{theorem}
\label{LW}\textbf{Part I : }Assume that $X$ is Gaussian, $\lambda\left(
\cdot\right)  $ is a convex decreasing function with $\lambda\left(  x\right)
\succ_{\infty}\exp\left(  -x^{\alpha}\right)  $ for some $\alpha>0$ and that
\ref{hyp-dens-eps} holds. Denote $T_{n}$ any estimator of the regression
function at a fixed point $r\left(  x_{0}\right)  =E\left(  y|X=x_{0}\right)
$ and $\mathcal{R}_{n}$ the minimax risk over the class $\mathcal{E}_{p}$
defined in (\ref{reg-class}) :%
\[
\mathcal{R}_{n}=\min_{T_{n}}\sup_{r\in\mathcal{E}_{p}}\mathbb{E}\left[
T_{n}-r\left(  x_{0}\right)  \right]  ^{2}%
\]
then $\mathcal{R}_{n}\succ\exp\left[  -c\left(  \log n\right)  ^{1-1/\left(
\alpha+1\right)  }\right]  $ which imples $n^{\beta}\mathcal{R}_{n}%
\rightarrow+\infty$ for any $\beta>0$ but $\left(  \log n\right)  ^{\beta
}\mathcal{R}_{n}\rightarrow0$ for any $\beta>0.$ Strengthening the assumptions
on $\lambda$ and taking $\lambda\left(  x\right)  \succ_{\infty}x^{-\alpha}$
for some $\alpha>1$ then $\mathcal{R}_{n}\succ\left(  \log n\right)
^{-\left(  \alpha-1\right)  }$.\newline\textbf{Part II : }Let\textbf{ }$X$ be
non gaussian but satisfy the conditions (\ref{A0}). Let $\rho$ be the
auxiliary function of the sall ball probability of $X$. Assume that $\rho$ is
regularly varying at $0$ with index $\alpha\geq1$ with either $\alpha>1$ or
$\alpha=1$ and $\rho\left(  s\right)  /s\succeq\log\left(  1/s\right)  $ then
again $n^{\beta}\mathcal{R}_{n}\rightarrow+\infty$ for any $\beta>0$.
\end{theorem}

In Part II we recall for the sake of completeness the conditions
$\mathbf{RV}_{\cdot}$ introduced sooner. The theorem above shows that it is
not possible to estimate the regression function in a nonparametric model with
functional inputs at a polynomial rate. The rates may be considered as
degenerate even when the functional variable $X$ is very smooth (case
$\lambda\left(  x\right)  =\exp\left(  -x^{\alpha}\right)  $ for some
$\alpha>0$) and the data concentrated close to a finite-dimensional space. In
the classical situations of polynomial decay, $\lambda\left(  x\right)  \simeq
x^{-\alpha}$ for some $\alpha>1$ the situation gets even worse and the optimal
rate we may recover is logarithmic. These negative results are clearly
connected with the complexity of the setting : the general nonparametric model
coupled with the sparsity of functional spaces already mentioned in the
paragraph below Proposition \ref{skye}.

\begin{remark}
Other classes of regression functions could be considered. Here $\mathcal{E}%
_{p}$ was considered because calculations are possible when looking for an
upper bound. However the theorem above holds, up to a change of constants when
$r$ blongs to a class $\mathcal{E}_{p}$ for which :%
\[
\sup_{r\in\mathcal{E}_{p}}\mathbb{E}\left(  r\left(  X\right)  -r\left(
x_{0}\right)  \right)  K\left(  \left\Vert X-x_{0}\right\Vert /h\right)
\asymp h^{p}F\left(  h\right)
\]
Like in a finite-dimensional framework, obtaining large values of $p$ switches
the problem to defining higher order kernels designed for functional data.
This issue is out of the scope of this work. Yet, because of the degeneracy of
the convergence rate we are not sure it deserves much attention in this setting.
\end{remark}

\textbf{Proof of Theorem \ref{LW}:}

The proof comes down to adapting Tsybakov (2004, Chapter 2, p.81) to our
framework. We consider two distant hypotheses : $r_{0}\equiv0$ and
$r_{1}\left(  x\right)  =2\left(  h^{\ast}\right)  ^{p}\mathcal{K}\left(
\left\Vert x-x_{0}\right\Vert /h\right)  $ with $\mathcal{K}\in\mathcal{E}%
_{p}$ and compactly supported. Here $h^{\ast}$ is the solution of the equation
(\ref{opt-band}). It is plain that $\left\vert r_{0}\left(  x_{0}\right)
-r_{1}\left(  x_{0}\right)  \right\vert =2\left(  h^{\ast}\right)  ^{p}$. Set
$z_{i}^{0}=\left(  y_{i}^{0},X_{i}\right)  $ and $z_{i}^{1}=\left(  y_{i}%
^{1},X_{i}\right)  $, denote $\mathbb{P}_{0}$ (resp. $\mathbb{P}_{1}$) the
distribution of the vector $\left(  z_{1}^{0},...,z_{n}^{0}\right)  $ (resp
$\left(  z_{1}^{1},...,z_{n}^{1}\right)  $) when the regression function is
$r_{0}$ (resp. $r_{1}$) and $\mathbb{P}_{0,i}$ (resp. $\mathbb{P}_{1,i}$) the
distribution of the margin $z_{i}^{0}$ (resp. $z_{i}^{1}$). We are going to
prove that the Hellinger-distance between $\mathbb{P}_{0}$ and $\mathbb{P}%
_{1}$ $\mathbf{H}\left(  \mathbb{P}_{0},\mathbb{P}_{1}\right)  $ is less than
a given $\tau<+\infty$. Let $f$ stand for the density of $U=\left\Vert
X-x_{0}\right\Vert /h$. The function $f$ is nothing but the first order
derivative of the small ball probability $F.$ We first compute the Hellinger
distance between the margins of $\mathbb{P}_{0}$ and $\mathbb{P}_{1}$ by
conditioning with respect to $U$. Let $\theta_{1}\left(  U\right)  =2\left(
h^{\ast}\right)  ^{p}\mathcal{K}\left(  U\right)  $ :%
\begin{align*}
\mathbf{H}^{2}\left(  \mathbb{P}_{0,i},\mathbb{P}_{1,i}\right)   &  \equiv\int%
{\displaystyle\int}
\left[  p_{\varepsilon}^{1/2}\left(  t\right)  -p_{\varepsilon}^{1/2}\left(
t-\theta_{1}\left(  u\right)  \right)  \right]  ^{2}f\left(  u\right)  dtdu\\
&  \leq p_{\ast}\int\theta_{1}^{2}\left(  u\right)  f\left(  u\right)
du=4p_{\ast}\left(  h^{\ast}\right)  ^{2p}\mathbb{E}\mathcal{K}^{2}\left(
U\right)
\end{align*}
by Assumption (\ref{hyp-dens-eps}). For $n$ large enough and by (\ref{sq-kern}%
) we deduce that :%
\begin{align*}
\mathbf{H}^{2}\left(  \mathbb{P}_{0,i},\mathbb{P}_{1,i}\right)   &
\leq8p_{\ast}\left(  h^{\ast}\right)  ^{2p}F\left(  h^{\ast}\right)
\mathcal{K}^{2}\left(  1\right) \\
&  \leq\mathbf{c}^{\ast}/n
\end{align*}

where $\mathbf{c}^{\ast}$ is some constant and we used (\ref{opt-band})$.$ The
decomposition of Hellinger distance for product measures (see Tsybakov (2004)
p. 69). gives
\begin{align*}
\mathbf{H}^{2}\left(  \mathbb{P}_{0},\mathbb{P}_{1}\right)   &  =2\left(
1-\left(  1-\frac{\mathbf{H}^{2}\left(  \mathbb{P}_{0,1},\mathbb{P}%
_{1,1}\right)  }{2}\right)  ^{n}\right) \\
&  \leq2\left(  1-\left(  1-\frac{\mathbf{c}^{\ast}}{2n}\right)  ^{n}\right)
\leq2\left(  1-\exp\left(  \frac{\mathbf{c}^{\ast}}{4}\right)  \right)
\end{align*}
and $\mathbf{H}^{2}\left(  \mathbb{P}_{0},\mathbb{P}_{1}\right)  \leq\tau$
with $\tau=2\left(  1-\exp\left(  \mathbf{c}^{\ast}/4\right)  \right)  $ which
almost finishes the proof of the Theorem. The last sentence is proved with the
same techniques and in view of Proposition \ref{eq}.

\section{Complementary facts}

In this short section are collected results of secondary interest. They
complete however the precedings by underlining some facts about the
non-unicity and the limits of the representation obtained above. Indeed the
preceding theorems lead to the following question : is it possible to obtain a
one to one representation, in a general framework, of the small ball
probabilities of random elements in $l_{2}$ -characterized by the sequence
$\left(  \lambda_{i}\right)  _{i\in\mathbb{N}}$- by a function in $\Gamma_{0}%
$, depending solely on its auxiliary function $\rho$ ? The answer is negative
for at least two reasons. First it is plain that two series $S$ and
$S^{\prime}$ built from different sequences $\left(  \lambda_{i},Z_{i}\right)
_{i\in\mathbb{N}}$ may have equivalent (at $0$) small ball probabilities.
Second, imagine that we confine to Gaussian small ball probabilities and
consider again the r.h.s. of (\ref{Lifs}) denoted $F\in\Gamma_{0}$ with
auxiliary function $\rho$. It is plain to see that any function $\phi F$ where
$\phi\left(  x+t\rho\left(  x\right)  \right)  /\phi\left(  x\right)
\rightarrow1$ when $x\rightarrow0$ belongs to $\Gamma_{0}$ with exactly the
same auxiliary function $\rho$. Consequently even fixing the distribution of
the sequence $Z_{i}$ is not sufficient to obtain a one to one mapping between
small ball probabilities and the set $\Gamma_{0}$.

Indeed, pick a function $F_{0}$ in the class $\Gamma_{0}.$ This function is
essentially defined by its auxiliary $\rho_{0}\left(  \cdot\right)  $ and
Theorem \ref{repres} is not precise enough for us to identify it with a small
ball probability. This is due to the non-unicity of $\rho$ mentioned just
under (\ref{zazie}) by the words "up to asymptotic equivalence". If $\rho
_{1}\sim_{0}\rho_{0}$ $\lim_{s\downarrow0^{+}}F_{0}\left(  s+x\rho_{1}\left(
s\right)  \right)  /F_{0}\left(  s\right)  =\exp\left(  x\right)  $ as well.
But the local behaviour at $0$ of $F_{1}\left(  s\right)  =\exp\left\{
\eta\left(  s\right)  -\int_{s}^{1}1/\rho_{1}\left(  t\right)  dt\right\}  $
may differ from $F_{0}\left(  s\right)  $ and $F_{0}$ may not be equivalent
with $F_{1}.$ What we show below is that if $F_{0}$ is accurately scaled we
may deduce from $F_{0}$ a new function $F_{0}^{\ast}$ which has the same
auxiliary function as $F_{0}$ (but which may not be equivalent to $F_{0}$) and
such that for a well-chosen sequence $\left(  \lambda_{i}\right)
_{i\in\mathbb{N}}$ and the Gaussian small ball probability $\mathbb{P}\left(
S<r\right)  $ is such that $\mathbb{P}\left(  S<r\right)  \sim_{0}F_{0}^{\ast
}\left(  r\right)  $

We start with a definition which seems to be new.

\begin{definition}
Let $\rho$ be a self-neglecting function. A function $\phi$ is called $\rho
$-self-neglecting if :%
\[
\frac{\phi\left(  x+t\rho\left(  x\right)  \right)  }{\phi\left(  x\right)
}\underset{x\rightarrow0}{\rightarrow}1.
\]

\end{definition}

It is obvious that, if $\phi$ is $\rho$-self-neglecting it is $\rho^{\ast}%
$-self-neglecting whenever $\rho^{\ast}\sim_{0}\rho$. We propose below in
Theorem (\ref{repres2}) a representation theorem for $\rho$-self-neglecting functions.

\begin{definition}
\label{eq-class}Pick a $\rho_{0}$ in the class of self-neglecting functions at
$0$ such that $\rho_{0}\left(  0\right)  =0$. We define the equivalence class
of a function $F_{0}\in\Gamma_{0}$ with auxiliary function $\rho_{0}$ by the
relationship $\triangle$ defined for all $G$ in $\Gamma_{0}$ by :%
\[
F_{0}\triangle G\Leftrightarrow\frac{F_{0}}{G}\mathrm{\ is\ }\rho
-\text{\textrm{self-neglecting for some }}\rho\sim_{0}\rho_{0}%
\]

\end{definition}

Remind that $\varphi\left(  t\right)  =t\gamma\left(  t\right)  $.

\begin{theorem}
\label{corsica}Let $F_{0}\in\Gamma_{0}$ with auxiliary function $\rho
_{0}=1/\gamma_{0}$. Assume that $\rho_{0}$ is regularly varying at $0$ with
index $\kappa$ $>1$ and $\mathbf{C}^{1}$ in a neighborhood of $0$. Consider
the equivalence class of $F_{0}$ in $\Gamma_{0}\backslash\triangle$ say
$\mathbf{F}_{0}$. Then one may pick $F_{0}^{\ast}\in\mathbf{F}_{0}$ such that
$F_{0}^{\ast}\left(  \cdot\right)  \sim_{0}\mathbb{P}\left(  S<\cdot\right)  $
were $S=\sum\lambda_{i}Z_{i}$, the $Z_{i}$'s follow a $\chi^{2}\left(
1\right)  $ distribution and :%
\[
\lambda_{i}=\frac{1}{\gamma\left(  \varphi^{-1}\left(  i\right)  \right)
}=\rho\left(  \varphi^{-1}\left(  i\right)  \right)
\]

\end{theorem}

\begin{remark}
Once again we encounter a regularly-varying condition on $\rho$. Here it
echoes in a way the assumption $\mathbf{A}_{0}$ (necessary to derive
(\ref{Dem})) which claims that the cdf of $Z$ is itself regularly varying at
$0$. An interesting open question would consist in finding examples of
auxiliary functions which are not regularly varying with positive index,
whenever it is possible.
\end{remark}

For the sake of completeness we obtain a last result, complementing and
illustrating Proposition \ref{eq-class}. From this Proposition we see that
$F\triangle G$ if $F=\phi G$ where $\phi$ is $\rho$-self-neglecting. The
forthcoming Theorem represents these functions $\phi$.

\begin{theorem}
\label{repres2}Let $\rho$ be self-neglecting at $0$ which does not vanish in a
neighborhood of $0$. A function $\phi$ is $\rho$-self-neglecting iff :%
\[
\phi\left(  x\right)  =c\left(  x\right)  \exp\left(  \int_{x}^{1}%
\frac{\varepsilon\left(  u\right)  }{\rho\left(  u\right)  }du\right)
\]
where $c\left(  u\right)  \rightarrow c\in\left]  0,+\infty\right)  $ and
$\varepsilon\left(  u\right)  \rightarrow0$ when $u\rightarrow0$ and
$\varepsilon$ has the same regularity as $\rho$.
\end{theorem}

This theorem generalizes the representation Theorem 2.11.3 for self-neglecting
functions p.121 in Bingham et al. (1987) initially due to Bloom (1976). If one
take $\phi=\rho$ the representation above coincides with the one announced in
this theorem.

\subsection{Conclusion and perspectives}

The first main results of this article identifies small ball probabilities in
$l_{2}$ with a class of rapidly varying functions involved in extreme value
theory and whose derivatives at all orders vanish at zero. This representation
was obtained through previous works especially the seminal and precious
formula (\ref{Dem}) of Lifshits (1997). We hope that this new formulation will
be more convenient for modelizing the small ball probabilities with some
applied -especially statistical- purposes in mind. However many other
questions arise. For instance the generalization to random elements with
values in $l_{p}$ or in more general Banach spaces is certainly an intricate
matter since the starting fomulas (\ref{Dem}) and followings seem to be
intimately suited to the space $l_{2}$.

A more promising track could be to explore the links between the auxiliary
function $\rho$, which inherits all the information on the regularity of $X$,
with the metric entropy of the unit ball of the reproducing kernel Hilbert
space of $X$ as explored in Li, Linde (1999) or with the degree of compactness
of the operator $v$ in Li, Linde (2004) for instance, the latter operator $v$
being obviously close to the covariance operator of $X$ hence in connection
with the $a_{i}$'s (or $\lambda_{i}$'s) of this article.

A surprising fact is the parallel that can be drawn between large deviations
on a one hand and extreme value theory on the other hand. Both were intially
introduced to model and explore large values of sequences of random elements.
It turns out that both provide an accurate setting to study small deviations
as well : Laplace transform for the classical approach and methods around the
domain of attraction of the third type (Gamma class, self-negclecting
functions...) as outlined here. However the connections between regular
variations and small ball probabilities have been known since de Bruin in
1959, and his theorem on Laplace transfoms (see Theorem 4.12.9 in Bingham et
al. (1987)). This work confirms that both Tauberian and extreme value theory
may provide tools complementing large deviations techniques to derive new
results in this area.

The other result shows, as an application of the previous, that the optimal
risk in nonparametric regression for functional data is degenerate in the
sense that we cannot expect obtain polynomial rates in the reasonable setting
used in this work. It is obviously interpretable in terms of curse of
dimensionality. A work is in progress to study the additive regression namely
the model :%
\[
y=\sum_{i=1}^{k}r_{i}\left(  \left\langle X,e_{i}\right\rangle \right)
+\varepsilon
\]
where the $r_{i}$ are functions defined on $\mathbb{R}$ and estimated from
one-dimensional projections of the data $X$. It is known since Stone(1985)
that this model is not subject to the curse of dimensionality when $X$ is
valued in $\mathbb{R}^{d}$. It would be a possible track to introduce
non-linearity in regression models for functional data and avoiding some
redhibitory features of a general model. The role of the auxiliary function
$\rho$ is major. The question of its estimation is quite simple indeed. From
Bingham et al. (1987) Corollary 3.10.5(b) p.177 we know that $\rho$ may be
taken as $F/F^{\prime}.$ A natural estimator of $\widehat{\rho}$ may be
$\widehat{F}/\widehat{f}$ where $\widehat{f}$ (resp. $\widehat{F}$) is a
kernel estimator of the density (resp. of the cumulative distribution
function) of $\left\Vert X\right\Vert .$ This is a simple procedure to check
some of the needed properties of $\rho$ such as its rate of decrease to $0.$

\section{Proofs}

Considerations about the smoothness at $0$ of $F$ and $\rho$ are not the
matter in this work and we will take it for granted that both functions are
smooth enough. Besides along the proofs we may sometimes consider generalized
or local inverses of some fonctions which may not be invertible or have smooth
derivatives everywhere. For example the auxiliary function $\rho$ defined on
$\mathbb{R}^{+}$ for which we always have $\rho^{\prime}\left(  0\right)  =0$
has no inverse on $\left[  0,c\right]  $ for $c>0$. But we may frequently use
the smoothness of, say, $\rho$ and $\rho^{-1}$ on sets $\left]  a,b\right[  $
for $0<a<b$ without always justifying it. We start with the proof of
Proposition \ref{cabrel}.

\textbf{Proof of Proposition \ref{cabrel} :}

Suppose that $\rho\left(  s\right)  /s$ does not tend to zero when $s$ does.
Then we may pick an $\varepsilon>0$ such that for infinitely many
$s_{k}\downarrow0$ when $k\uparrow+\infty,$ $\rho\left(  s_{k}\right)
/s_{k}>\varepsilon$. Now fix $x<-\varepsilon^{-1}$ then $s_{k}+x\rho\left(
s_{k}\right)  <0$ and $F\left(  s_{k}+x\rho\left(  s_{k}\right)  \right)  =0$
for all $k$ and $F\left(  s_{k}+x\rho\left(  s_{k}\right)  \right)  /F\left(
s_{k}\right)  $ cannot converge to $\exp\left(  x\right)  $. The second part
of the proof, namely ensuring the $\rho$ is self-neglecting, follows the lines
of the proof of Proposition 3.10.6 in Bingham et al. (1987).\endproof\medskip

\textbf{Proof of Proposition \ref{skye}} : Suppose that for some $p$
$F^{\left(  p\right)  }\left(  0\right)  \neq0$ and take $p^{\ast}%
=\inf\left\{  p\in\mathbb{N}:F^{\left(  p\right)  }\left(  0\right)
\neq0\right\}  $. It is plain that $F^{\left(  p^{\ast}\right)  }\left(
0\right)  >0$ since $F$ is positive. Then we should consider two cases. First
if $F^{\left(  p^{\ast}\right)  }\left(  0\right)  =c<+\infty$ then $F\left(
s\right)  \sim cs^{p^{\ast}}$. Taking :%
\[
\frac{F\left(  s+\rho\left(  s\right)  \right)  }{F\left(  s\right)  }%
=\frac{F\left(  s+\rho\left(  s\right)  \right)  }{\left(  s+\rho\left(
s\right)  \right)  ^{p^{\ast}}}\frac{s^{p^{\ast}}}{F\left(  s\right)  }%
\frac{\left(  s+\rho\left(  s\right)  \right)  ^{p^{\ast}}}{s^{p^{\ast}}}%
\]
we see that the left hand side of the display above tends to $\exp\left(
1\right)  $ whereas the right hand side tends to $1$.

Second if $F^{\left(  p^{\ast}\right)  }\left(  0\right)  =+\infty$ we clearly
have $F\left(  s\right)  /s^{p^{\ast}}\rightarrow+\infty$ when $s\rightarrow
0$. Take $\varepsilon$ such that $1/\varepsilon>p^{\ast}+2$. Since
$\rho^{\prime}\left(  0\right)  =0$ and $\rho$ is positive we may pick an
$s_{0}$ such that $\sup_{0\leq u\leq s_{0}}\rho^{\prime}\left(  u\right)
\leq\varepsilon$. From (\ref{zazie}) we get :%
\begin{align*}
\frac{F\left(  s\right)  }{s^{p}}  &  \leq\frac{C}{s^{p}\rho^{2}\left(
s\right)  }\exp\left\{  -\int_{s}^{1}\frac{1}{\rho\left(  t\right)
}dt\right\}  \leq\frac{C}{\rho^{2+p}\left(  s\right)  }\exp\left\{  -\int
_{s}^{1}\frac{1}{\rho\left(  t\right)  }dt\right\} \\
&  \leq\frac{C^{\prime}}{\rho^{2+p}\left(  s\right)  }\exp\left\{  -\int
_{s}^{s_{0}}\frac{1}{\rho\left(  t\right)  }dt\right\}
\end{align*}
where we assume that $s\leq s_{0}$. Then we have%
\begin{align*}
\exp\left\{  -\int_{s}^{s_{0}}\frac{1}{\rho\left(  t\right)  }dt\right\}   &
=\exp\left\{  -\int_{s}^{s_{0}}\frac{\rho^{\prime}\left(  t\right)  }%
{\rho\left(  t\right)  }\frac{1}{\rho^{\prime}\left(  t\right)  }dt\right\}
\leq\exp\left\{  -\frac{1}{\varepsilon}\int_{s}^{s_{0}}\frac{\rho^{\prime
}\left(  t\right)  }{\rho\left(  t\right)  }dt\right\} \\
&  =\exp\left\{  \frac{1}{\varepsilon}\ln\rho\left(  s\right)  -\frac
{1}{\varepsilon}\ln\rho\left(  s_{0}\right)  \right\}  .
\end{align*}
At last%
\[
\frac{F\left(  s\right)  }{s^{p}}\leq C^{\prime\prime}\left[  \rho\left(
s\right)  \right]  ^{\frac{1}{\varepsilon}-p^{\ast}-2}%
\]
which contradicts the fact that $F\left(  s\right)  /s^{p^{\ast}}%
\rightarrow+\infty$.\endproof\medskip

We start the proof of Theorem \ref{main}\textbf{ }

\textbf{Proof of Theorem \ref{main} :}

From Definition \ref{def} and (\ref{Dem}) we see that Theorem \ref{main} holds
whenever for all $x\in\mathbb{R}$ :%
\[
\lim_{s\rightarrow0}\frac{\gamma\left(  s\right)  \sigma\left(  s\right)
}{\gamma\left(  s+x\rho\left(  s\right)  \right)  \sigma\left(  s+x\rho\left(
s\right)  \right)  }\exp\left(  \left(  s+x\rho\left(  s\right)  \right)
\gamma\left(  s+x\rho\left(  s\right)  \right)  -s\gamma\left(  s\right)
\right)  \frac{\Lambda\left(  \gamma\left(  s+x\rho\left(  s\right)  \right)
\right)  }{\Lambda\left(  \gamma\left(  s\right)  \right)  }=\exp x.
\]
We will more specifically prove below that when $s$ decays to $0$ :%
\begin{align*}
\gamma\left(  s+x\rho\left(  s\right)  \right)  /\gamma\left(  s\right)
\sigma\left(  s+x\rho\left(  s\right)  \right)  /\sigma\left(  s\right)   &
\rightarrow1\\
\exp\left(  \left(  s+x\rho\left(  s\right)  \right)  \gamma\left(
s+x\rho\left(  s\right)  \right)  -s\gamma\left(  s\right)  -x\right)
\frac{\Lambda\left(  \gamma\left(  s+x\rho\left(  s\right)  \right)  \right)
}{\Lambda\left(  \gamma\left(  s\right)  \right)  }  &  \rightarrow1
\end{align*}
The two next lemmas are dedicated to showing that, in the above display the
fraction as well as the exponential both tend to $1$ when $s$ goes to zero and
$\rho$ is chosen as in the Theorem.We just have to clarifiy formula
(\ref{Lifs}) within the Theorem. This stem directly from (\ref{Dem}). Indeed
from (\ref{mu}) and (\ref{D4}) we see that $\sigma^{2}=-\partial
r/\partial\gamma$ and we just have to show that $\gamma r+\log\Lambda\left(
\gamma\right)  =\int_{r_{0}}^{r}\gamma\left(  s\right)  ds$. Elementary
calculations yield :%
\[
\frac{\partial\left(  \gamma r+\log\Lambda\left(  \gamma\right)  \right)
}{\partial r}=\gamma\left(  r\right)  .
\]
Let $r_{0}=\mathbb{E}Z\cdot\sum_{j=1}^{n}\lambda_{j}.$ Applying formula
(\ref{mu}) at $\gamma=0$ we notice that $\gamma\left(  r_{0}\right)
=0=\log\Lambda\left(  \gamma\left(  r_{0}\right)  \right)  $ and we
conclude.\medskip

\begin{lemma}
\label{L1}Take $\rho\left(  s\right)  =1/\gamma\left(  s\right)  ,$ then :%
\[
\lim_{s\rightarrow0}\exp\left(  \left(  s+x\rho\left(  s\right)  \right)
\gamma\left(  s+x\rho\left(  s\right)  \right)  -s\gamma\left(  s\right)
-x\right)  \frac{\Lambda\left(  \gamma\left(  s+x\rho\left(  s\right)
\right)  \right)  }{\Lambda\left(  \gamma\left(  s\right)  \right)  }=1.
\]
\bigskip
\end{lemma}

\begin{remark}
\label{samba}Obviously $\gamma$ has at least two (we do not need more)
continuous derivatives on a neighborhood of infinity (here $\left]
1,+\infty\right)  $ for instance). It is also strightforward to see that
$\gamma$, which is strictly decreasing on $\left]  1,+\infty\right)  ,$ is
also a $C^{1}$ diffeomorphism on this set. Clearly $\lim_{s\rightarrow0}%
\rho\left(  s\right)  =0$ but from (\ref{Boine}) it is plain that $\rho\left(
s\right)  /s$ also tends to zero when $s$ does which implies that
$\rho^{\prime}\left(  0\right)  =0.$ Indeed proving that $\rho\left(
s\right)  /s$ tends to zero comes down to proving that $s\gamma\left(
s\right)  \rightarrow+\infty$.\medskip
\end{remark}

\textbf{Proof of Lemma \ref{L1} :}

Denote $I\left(  s\right)  =s\gamma\left(  s\right)  +\log\Lambda\left(
\gamma\left(  s\right)  \right)  $. We should prove that :%
\[
\lim_{s\rightarrow0}I\left(  s+x\rho\left(  s\right)  \right)  -I\left(
s\right)  -x=0
\]
Taylor's formula gives :%
\begin{equation}
I\left(  s+x\rho\left(  s\right)  \right)  -I\left(  s\right)  =x\rho\left(
s\right)  I^{\prime}\left(  s\right)  +\frac{x^{2}}{2}\rho^{2}\left(
s\right)  I^{\prime\prime}\left(  c_{s,x}\right)  \label{signs}%
\end{equation}
where $c_{s,x}=c$ lies somewhere in $\left[  s,s+x\rho\left(  s\right)
\right]  $ if $x\geq0$ and in $\left[  s+x\rho\left(  s\right)  ,x\right]  $
if $x<0$. From (\ref{D4}) we see that :%
\begin{align*}
I^{\prime}\left(  s\right)   &  =\gamma\left(  s\right)  +s\gamma^{\prime
}\left(  s\right)  +\gamma^{\prime}\left(  s\right)  \frac{\partial\log
\Lambda\left(  \gamma\left(  s\right)  \right)  }{\partial\gamma}\\
&  =\gamma\left(  s\right)
\end{align*}
Hence (\ref{signs}) may be rewritten :%
\begin{align*}
I\left(  s+x\rho\left(  s\right)  \right)  -I\left(  s\right)   &
=x+\frac{x^{2}}{2}\rho^{2}\left(  s\right)  \gamma^{\prime}\left(
c_{s,x}\right) \\
&  =x+\frac{x^{2}}{2}\frac{\gamma^{\prime}\left(  c_{s,x}\right)  }{\gamma
^{2}\left(  s\right)  }=x-\frac{x^{2}}{2}\frac{d\left(  1/\gamma\right)  }%
{ds}\left(  c_{s,x}\right)  \cdot\frac{\gamma^{2}\left(  c_{s,x}\right)
}{\gamma^{2}\left(  s\right)  }%
\end{align*}

We first show that $\gamma^{2}\left(  c_{s,x}\right)  /\gamma^{2}\left(
s\right)  =\rho^{2}\left(  s\right)  /\rho^{2}\left(  c_{s,x}\right)  $ is
bounded above. We may always write $c_{s,x}=s+t_{x}\left(  s\right)
\rho\left(  s\right)  $ where $-x\leq t_{x}\left(  s\right)  \leq x$ for all
$s$. Taylor's formula yields%
\[
\rho\left(  s+t_{x}\left(  s\right)  \rho\left(  s\right)  \right)
=\rho\left(  s\right)  +t_{x}\left(  s\right)  \rho\left(  s\right)
\rho^{\prime}\left(  d\right)  =\rho\left(  s\right)  \left(  1+t_{x}\left(
s\right)  \rho^{\prime}\left(  d\right)  \right)
\]
where $d$ lies between $s$ and $s+t_{x}\left(  s\right)  \rho\left(  s\right)
$. Hence :%
\[
\frac{\rho\left(  s\right)  }{\rho\left(  c_{s,x}\right)  }=\frac{1}%
{1+t_{x}\left(  s\right)  \rho^{\prime}\left(  d\right)  }\leq\frac
{1}{1-\left\vert x\right\vert \rho^{\prime}\left(  d\right)  }%
\]
The continuity of $\rho^{\prime}$ at $0$ and its nullity at $0$ (see Remark
\ref{samba}) implies on a one hand that the display above is bounded above for
fixed $x$ and $s$ (hence $d$) going to zero and also that :%
\[
\frac{d\left(  1/\gamma\right)  }{ds}\left(  c_{s,x}\right)  =\rho^{\prime
}\left(  c_{s,x}\right)  \rightarrow0.
\]
At last, $I\left(  s+x\rho\left(  s\right)  \right)  -I\left(  s\right)
\rightarrow x$ which finishes the proof of the Lemma.\endproof.\bigskip

\begin{lemma}
\label{AboutHer}We have :%
\[
\lim_{s\rightarrow0}\frac{\gamma\left(  s+x\rho\left(  s\right)  \right)
\sigma\left(  s+x\rho\left(  s\right)  \right)  }{\gamma\left(  s\right)
\sigma\left(  s\right)  }=1
\]

\end{lemma}

\textbf{Proof :}

Once again Taylor's formula leads to :%
\begin{equation}
\frac{\gamma\left(  s+x\rho\left(  s\right)  \right)  \sigma\left(
s+x\rho\left(  s\right)  \right)  -\gamma\left(  s\right)  \sigma\left(
s\right)  }{\gamma\left(  s\right)  \sigma\left(  s\right)  }=\frac
{x\rho\left(  s\right)  }{\gamma\left(  s\right)  \sigma\left(  s\right)
}\left[  \gamma^{\prime}\left(  c\right)  \sigma\left(  c\right)
+\gamma\left(  c\right)  \sigma^{\prime}\left(  c\right)  \right]
\label{brass}%
\end{equation}

where $c\in\left(  s,s\pm x\rho\left(  s\right)  \right)  $. We will prove
that $\left[  \gamma^{\prime}\left(  c\right)  \sigma\left(  c\right)
+\gamma\left(  c\right)  \sigma^{\prime}\left(  c\right)  \right]  /\gamma
^{2}\left(  s\right)  \sigma\left(  s\right)  $ tends to zero. We cut the
latter into two terms. First consider%
\[
\frac{\rho\left(  s\right)  }{\gamma\left(  s\right)  \sigma\left(  s\right)
}\gamma^{\prime}\left(  c\right)  \sigma\left(  c\right)  =\frac
{\gamma^{\prime}\left(  c\right)  }{\gamma^{2}\left(  c\right)  }\frac
{\gamma^{2}\left(  c\right)  }{\gamma^{2}\left(  s\right)  }\frac
{\sigma\left(  c\right)  }{\sigma\left(  s\right)  }%
\]
We proved above within the proof of the previous Lemma (\ref{L1}) that
$\gamma^{2}\left(  c\right)  /\gamma^{2}\left(  s\right)  $ is bounded above.
We proved as well that $\gamma^{\prime}\left(  c\right)  /\gamma^{2}\left(
c\right)  $ tends to zero when $c$ does. Finally we should just control
$\sigma\left(  c\right)  /\sigma\left(  s\right)  .$ We have $\sigma\left(
c\right)  =\sigma\left(  s\right)  +\left(  c-s\right)  \sigma^{\prime}\left(
\xi\right)  $ where $\xi\in\left[  s,c\right]  $ hence%
\[
0\leq\frac{\sigma\left(  c\right)  }{\sigma\left(  s\right)  }=1+\frac
{c-s}{\sigma\left(  s\right)  }\sigma^{\prime}\left(  \xi\right)  \leq
1+\frac{x}{\gamma\left(  s\right)  \sigma\left(  s\right)  }\sigma^{\prime
}\left(  \xi\right)  .
\]
We see in Lifshits (1997, Lemma 2 p.431) that $\lim_{s\rightarrow0}%
\gamma\left(  s\right)  \sigma\left(  s\right)  =+\infty$ and that
$\sigma\left(  s\right)  \leq sc_{13}^{-1}$ where $c_{13}$ is some constant
from which it is plain that $\sup_{\xi\in\mathcal{V}_{0}}\left\vert
\sigma^{\prime}\left(  \xi\right)  \right\vert <+\infty$ where $\mathcal{V}%
_{0}$ is any neighborhood of $0$. We deduce that $\sigma^{\prime}\left(
\xi\right)  /\gamma\left(  s\right)  \sigma\left(  s\right)  $ tends to zero
which finally yields%
\[
\frac{\rho\left(  s\right)  }{\gamma\left(  s\right)  \sigma\left(  s\right)
}\gamma^{\prime}\left(  c\right)  \sigma\left(  c\right)  \rightarrow0.
\]
We turn to the second term in (\ref{brass}) : $\rho\left(  s\right)
\gamma\left(  c\right)  \sigma^{\prime}\left(  c\right)  /\gamma\left(
s\right)  \sigma\left(  s\right)  .$ We rewrite it :%
\[
\frac{\gamma\left(  c\right)  \sigma^{\prime}\left(  c\right)  }{\gamma
^{2}\left(  s\right)  \sigma\left(  s\right)  }=\frac{1}{\gamma\left(
s\right)  \sigma\left(  s\right)  }\frac{\gamma\left(  c\right)  }%
{\gamma\left(  s\right)  }\sigma^{\prime}\left(  c\right)
\]
As shown above from Lisfhits' work : $\gamma\sigma\rightarrow+\infty,$
$\sup_{c\in\mathcal{V}_{0}}\left\vert \sigma^{\prime}\left(  c\right)
\right\vert <+\infty$ and $\gamma\left(  c\right)  /\gamma\left(  s\right)  $
is bounded above and this second term also decays to zero. This finishes the
proof of Lemma \ref{AboutHer}.\bigskip Now we turn to the proof of the
converse part, Theorem \ref{corsica}. It takes two steps.

First we should make sure that when $\lambda_{i}=\rho\left(  \varphi
^{-1}\left(  i\right)  \right)  ,\sum\lambda_{i}<+\infty$ which will ensure
that the random element defined by $S=\sum\lambda_{i}Z_{i}$ is well-defined.

\begin{lemma}
\label{rouge}When $\lambda_{i}=\rho\left(  \varphi^{-1}\left(  i\right)
\right)  ,$ $\sum\lambda_{i}<+\infty$.\medskip
\end{lemma}

\textbf{Proof :} It is easily seen that $\varphi^{-1}$ is non decreasing in a
neighborhood of $+\infty$. Indeed it suffices to prove that $\varphi$ is,
which may be deduced from its definition by studying its derivative. By the
way one may also see that $\varphi$ is concave. Now since $\varphi^{-1}$ is
non decreasing it is enough to prove that :%
\[
\int^{+\infty}\rho\left(  \varphi^{-1}\left(  x\right)  \right)  dx<+\infty
\]
where the notation above means "the improper integral converges at infinity".
Set $u=\varphi^{-1}\left(  x\right)  $ above then we should examine :%
\[
\int_{0}\rho\left(  u\right)  \varphi^{\prime}\left(  u\right)  du.
\]
Integrating by part this comes down to ensuring first that $\rho\left(
u\right)  \varphi\left(  u\right)  =u$ tends to a finite limit as $u$ tends to
$0$ which is plain and that
\[
\int_{0}\rho^{\prime}\left(  u\right)  \varphi\left(  u\right)  du=\int
_{0}u\frac{\rho^{\prime}\left(  u\right)  }{\rho\left(  u\right)  }du<+\infty
\]
Now we are in a position to apply Karamata's theorem to $\rho^{\prime}$ :
since $\rho$ is regularly varying at $0$ with index $d\geq1$ (since
$\rho^{\prime}\left(  0\right)  =0$), and monotone in a right neighborhood of
zero, $\rho^{\prime}$ is also regularly varying with index $\geq0$ (see
Theorem 1.732.b p.39 in Bingham et al.(1987)). Then we can apply the direct
part of Karamata's Theorem to $\rho^{\prime}$ (see ibid. Theorem 1.5.11 (i)
p.28 where the limit should be taken here at zero) and%
\[
\lim_{t\rightarrow0}\frac{t\rho^{\prime}\left(  t\right)  }{\rho\left(
t\right)  }<+\infty
\]
which ensures that the integral above converges and finally that $\sum
\lambda_{i}<+\infty$. this completes the proof of Lemma \ref{rouge}%
.\endproof

\textbf{Proof of Theorem \ref{corsica} :}

Pick an $F_{0}$ in $\Gamma_{0}$ with auxiliary function $\rho_{0}$ and
consider the function $F_{0}^{\ast}\left(  r\right)  =\sqrt{\rho_{0}^{\prime
}\left(  r\right)  /\pi}\exp\left[  -\int_{r}^{r_{0}}ds/\rho_{0}\left(
s\right)  \right]  $ with $r_{0}=\sum_{i}\rho\left(  \varphi^{-1}\left(
i\right)  \right)  .$ Note that $\sqrt{\rho^{\prime}\left(  \cdot\right)  }$
hence $\rho^{\prime}\left(  \cdot\right)  $ are $\rho$-self-neglecting because
:%
\[
\frac{\rho^{\prime}\left(  r+x\rho\left(  r\right)  \right)  }{\rho^{\prime
}\left(  r\right)  }\rightarrow_{r\rightarrow0}1
\]
Indeed $\rho^{\prime}\left(  r+x\rho\left(  r\right)  \right)  =\rho^{\prime
}\left(  r\left(  1+x\rho\left(  r\right)  /r\right)  \right)  $,
$\rho^{\prime}$ is regularly varying with positive index since $\rho$ is
itself regularly varying with index $\kappa>1,$ and $\rho\left(  r\right)
/r\rightarrow0$ lead to
\[
\lim_{r\rightarrow0}\rho^{\prime}\left(  r\left(  1+x\rho\left(  r\right)
/r\right)  \right)  /\rho^{\prime}\left(  r\right)  =\lim\left(
1+x\rho\left(  r\right)  /r\right)  ^{\kappa-1}=1
\]
This proves that $F_{0}^{\ast}\vartriangle F_{0}$. It remains to show that
$F_{0}^{\ast}\sim_{0}\mathbb{P}\left(  S<\cdot\right)  $. Like above
$\gamma_{0}=1/\rho_{0}$. Start from (\ref{log-lap-gauss}) that is $r=\sum
_{j}\lambda_{j}/\left(  1+2\gamma_{0}\lambda_{j}\right)  .$ Now, following the
proof of Proposition \ref{eq} we set $J\left(  r\right)  =r/\rho_{0}\left(
r\right)  $ (we just make use of display (\ref{sc}), fix $J\left(  r\right)
\rho_{0}\left(  r\right)  /r=1$ instead of bounding it above and below) and
take $a\left(  \cdot\right)  =J^{-1}\left(  \cdot\right)  $ then finally
$S=\sum_{i=1}^{+\infty}Z_{i}/a\left(  i\right)  $. By construction
$\mathbb{P}\left(  S<\cdot\right)  \sim F_{0}^{\ast}$.

Finally we turn to the proof of Theorem \ref{repres2} and start with a Lemma.
This Lemma, its proof and the subsequent proof of the theorem adapt the
derivation of Lemma 2.11.2 and Theorem 2.11.3 of Bingham et al. (1987).

\begin{lemma}
Let $\rho$ be self-neglecting at $0.$ For $x_{0}>0$ sufficiently small the
sequence $x_{n}=x_{n-1}-\rho\left(  x_{n-1}\right)  $ tends to $0$.\medskip
\end{lemma}

\textbf{Proof :} First note that the sequence $x_{n}$ is decreasing since
$\rho\geq0$ and notice from the properties of self-neglecting functions
(namely $\rho\left(  s\right)  /s\rightarrow0$ when $s\rightarrow0$) that for
a sufficently small $x_{0}>0,$ $x_{n}\geq0$ for all $n$. The limit of $x_{n}$
exists, is denoted $l$. Suppose that $l>0$. Then $\rho\left(  l\right)  >0$
and since $\rho$ is a non decreasing function $\rho\left(  x_{k}\right)
\geq\rho\left(  l\right)  $ for all $k$. At last%
\begin{align*}
x_{n}  &  =x_{n-1}-\rho\left(  x_{n-1}\right)  =x_{0}-\sum_{k=0}^{n-1}%
\rho\left(  x_{k}\right) \\
&  \leq x_{0}-n\rho\left(  l\right)  .
\end{align*}
Letting $n$ go to infinity $x_{n}$ goes to $-\infty$ which contradicts
$x_{n}\geq0$ hence the Lemma.\endproof\medskip

\textbf{Proof of Theorem \ref{repres2}:}

Let $x_{n}$ be as in the preceding Lemma. Let $p$ be a $C^{\infty}$
probability density on $\left[  0,1\right]  $ and set for $x_{n+1}\leq u\leq
x_{n}$%
\[
\varepsilon\left(  u\right)  =\frac{\ln\phi\left(  x_{n+1}\right)  -\ln
\phi\left(  x_{n}\right)  }{x_{n}-x_{n+1}}p\left(  \frac{x_{n}-u}%
{x_{n}-x_{n+1}}\right)  \rho\left(  u\right)  .
\]

The proof takes three steps.

We prove first that for all $x_{n},$ $\phi\left(  x_{n}\right)  =\exp\left(
\int_{x_{n}}^{1}\frac{\varepsilon\left(  u\right)  }{\rho\left(  u\right)
}du\right)  .$ In fact we may always define $\varepsilon\left(  u\right)  $,
$x_{0}\leq u\leq1$ such that $\phi\left(  x_{0}\right)  =\exp\left(
\int_{x_{0}}^{1}\frac{\varepsilon\left(  u\right)  }{\rho\left(  u\right)
}du\right)  .$ Then assume that $\phi\left(  x_{k}\right)  =\exp\left(
\int_{x_{k}}^{1}\frac{\varepsilon\left(  u\right)  }{\rho\left(  u\right)
}du\right)  $ for $k=0,1,..,n$. We have :%
\begin{align*}
\int_{x_{n+1}}^{1}\frac{\varepsilon\left(  u\right)  }{\rho\left(  u\right)
}du  &  =\int_{x_{n+1}}^{x_{n}}\frac{\varepsilon\left(  u\right)  }%
{\rho\left(  u\right)  }du+\int_{x_{n}}^{1}\frac{\varepsilon\left(  u\right)
}{\rho\left(  u\right)  }du\\
&  =\ln\phi\left(  x_{n}\right)  +\frac{\ln\phi\left(  x_{n+1}\right)
-\ln\phi\left(  x_{n}\right)  }{x_{n}-x_{n+1}}\int_{x_{n+1}}^{x_{n}}p\left(
\frac{x_{n}-u}{x_{n}-x_{n+1}}\right)  du\\
&  =\ln\phi\left(  x_{n}\right)  -\left(  \ln\phi\left(  x_{n+1}\right)
-\ln\phi\left(  x_{n}\right)  \right)  \int_{1}^{0}p\left(  t\right)  dt\\
&  =\ln\phi\left(  x_{n+1}\right)
\end{align*}

Second we prove that for $x_{n+1}\leq x\leq x_{n}$ $\lim_{x\rightarrow0}%
\phi\left(  x\right)  /\phi\left(  x_{n}\right)  =1$. We note that
$x=x_{n}-\lambda_{x}\rho\left(  x_{n}\right)  $ where $\lambda_{x}\in\left[
0,1\right]  $ hence%
\[
\lim_{x\rightarrow0}\frac{\phi\left(  x_{n}-\lambda_{x}\rho\left(
x_{n}\right)  \right)  }{\phi\left(  x_{n}\right)  }=1
\]
uniformly with respect to $\lambda_{x}\in\left[  0,1\right]  $.

The third and last step is devoted to proving that $\left\vert \varepsilon
\left(  u\right)  \right\vert \rightarrow0$ when $u\rightarrow0.$ Indeed for
all $x_{n+1}\leq u\leq x_{n},$%
\[
\left\vert \varepsilon\left(  u\right)  \right\vert \leq\left\vert
p\right\vert _{\infty}\left\vert \frac{\ln\phi\left(  x_{n+1}\right)  -\ln
\phi\left(  x_{n}\right)  }{x_{n}-x_{n+1}}\right\vert \rho\left(  u\right)
\]
We focus on%
\begin{align*}
\left\vert \frac{\ln\phi\left(  x_{n+1}\right)  -\ln\phi\left(  x_{n}\right)
}{x_{n}-x_{n+1}}\right\vert \rho\left(  u\right)   &  =\frac{\rho\left(
u\right)  }{\rho\left(  x_{n}\right)  }\ln\frac{\phi\left(  x_{n}\right)
}{\phi\left(  x_{n+1}\right)  }\\
&  =\frac{\rho\left(  x_{n}-\lambda_{u}\rho\left(  x_{n}\right)  \right)
}{\rho\left(  x_{n}\right)  }\ln\frac{\phi\left(  x_{n}\right)  }{\phi\left(
x_{n}-\rho\left(  x_{n}\right)  \right)  }%
\end{align*}
Just like above $\rho\left(  x_{n}-\lambda_{u}\rho\left(  x_{n}\right)
\right)  /\rho\left(  x_{n}\right)  \rightarrow1$ since $\rho$ is
self-neglecting. Finally by the definition of $\phi$ we get%
\[
\ln\frac{\phi\left(  x_{n}\right)  }{\phi\left(  x_{n}-\rho\left(
x_{n}\right)  \right)  }\rightarrow0
\]
which finishes the proof of the Theorem.\endproof

\end{document}